\documentclass[format=acmsmall]{acmart}

\usepackage{amsmath}
\usepackage{amsfonts}
\usepackage{graphicx}
\usepackage{amsthm}
\usepackage{amscd}
\usepackage{amsmath}
\usepackage{latexsym}
\usepackage{xparse}
\usepackage{xifthen}
\usepackage{subfigure}
\usepackage[scaled=.8]{beramono}
\usepackage[boxed, lined]{algorithm2e}
\usepackage{tikz-uml}

\acmJournal{TOMS}
\acmVolume{V}
\acmNumber{N}
\acmArticle{A}
\acmYear{YYYY}
\acmMonth{1}
\copyrightyear{YYYY}

\setcopyright{acmcopyright}

\acmDOI{0000001.0000001}

\SetKwFor{ForM}{for}{}{end}
\SetKwProg{Fn}{function}{}{end}
\SetKwInOut{In}{input}
\SetKwInOut{Out}{output}

\newcommand{\mbf}[1]{{\boldsymbol{#1}}}
\newcommand{\supp}{\mathrm{supp}}

\newcommand{\RR}{{\mathbb{R}}}
\newcommand{\ZZ}{{\mathbb{Z}}}
\newcommand{\CC}{{\mathbb{C}}}

\newcommand{\NN}{{\mathbb{N}}}
\newcommand{\GG}{{\mathbb{G}}}
\newcommand{\PP}{{\mathbb{P}}}
\newcommand{\EE}{{\mathbb{GE}}}
\newcommand{\TT}{{\mathbb{GT}}}
\newcommand{\PE}{{\mathbb{PE}}}
\newcommand{\PT}{{\mathbb{PT}}}
\newcommand{\WW}{{\mathbb{W}}}

\newcommand{\domain}{\Delta}
\newcommand{\interval}{J}
\newcommand{\nelms}{m}
\newcommand{\nsum}{\mu}
\newcommand{\croot}{\omega}

\renewcommand{\SS}[2]{{\mathbb{S}_{#1}^{#2}}}

\newcommand{\ECT}[2]{{\mathbb{T}_{#1}^{#2}}}
\newcommand{\uu}{g}
\newcommand{\bs}{B}
\newcommand{\BS}{N}
\newcommand{\ts}{\phi}

\renewcommand{\r}{r}
\newcommand{\br}{\mbf{r}}
\newcommand{\deriv}{k}

\newcommand{\ii}{{\mathrm{i}}}
\newcommand{\ee}{{\mathrm{e}}}

\renewcommand{\a}{a}
\renewcommand{\b}{b}
\renewcommand{\aa}{z}
\newcommand{\dd}{{\,\mathrm{d}}}

\newcommand{\x}{x}

\newcommand{\y}{y}
\newcommand{\dy}{\dd y}

\newcommand{\N}{n}
\newcommand{\p}{p}
\newcommand{\q}{q}
\newcommand{\bp}{\mbf{\p}}
\newcommand{\knot}{\xi}
\newcommand{\lknot}{u}
\newcommand{\rknot}{v}
\newcommand{\bspint}{d}
\newcommand{\bbint}{b}
\newcommand{\lintsum}{\mu_{\mbf{\rknot}}}
\newcommand{\rintsum}{\mu_{\mbf{\lknot}}}
\newcommand{\ca}{\alpha}
\newcommand{\cb}{\beta}
\newcommand{\cL}{\mathcal{L}}

\newcommand{\splSpacep}{\SS{\bp}{}(\domain)}
\newcommand{\splSpacerp}{\SS{\bp}{\br}(\domain)}
\newcommand{\splSpacerptilde}{\SS{\bp}{\tilde{\br}}(\domain)}

\newcommand{\Chef}{Tchebycheffian}
\newcommand{\Chefshort}{Tchebycheff}

\newcommand{\Matlab}{\textsc{Matlab}}

\newtheorem{proposition}{Proposition}
\newtheorem{remark}{Remark}
\newtheorem{example}{Example}
\newtheorem{definition}{Definition}

\begin{document}
\title{Algorithm xxx: Computation of Multi-Degree \Chef{} B-Splines}  

\author{Hendrik Speleers}
\affiliation{%
  \institution{Department of Mathematics, University of Rome Tor Vergata}
  \streetaddress{Via della Ricerca Scientifica 1}
  \city{Rome}
  \postcode{00133}
  \country{Italy}}
\email{speleers@mat.uniroma2.it}

\begin{abstract}
Multi-degree \Chef{} splines are splines with pieces drawn from extended (complete) \Chefshort{} spaces, which may differ from interval to interval, and possibly of different dimensions. These are a natural extension of multi-degree polynomial splines. Under quite mild assumptions, they can be represented in terms of a so-called MDTB-spline basis; such basis possesses all the characterizing properties of the classical polynomial B-spline basis. We present a practical framework to compute MDTB-splines, and provide an object-oriented implementation in \Matlab{}. The implementation supports the construction, differentiation, and visualization of MDTB-splines whose pieces belong to \Chefshort{} spaces that are null-spaces of constant-coefficient linear differential operators. The construction relies on an extraction operator that maps local \Chef{} Bernstein functions to the MDTB-spline basis of interest.
\end{abstract}

%
%

%


\maketitle


\section{Introduction}

Splines are undoubted an important tool in several branches of the sciences including geometric modeling, signal processing, data analysis, visualization, and numerical simulation, just to mention a few \cite{Cohen:2001,Cottrell:2009}.
The term \emph{splines} usually refers to univariate piecewise (algebraic) polynomial functions with certain smoothness, whose popularity can be mainly attributed to their representation in terms of the so-called \emph{B-splines}. The B-splines enjoy properties as local linear independence, minimal support, non-negativity and partition of unity; they can be computed through a stable recurrence relation; and they can even be seen as the geometrically optimal basis for piecewise polynomial spaces.

This raises the following natural question: is there a more general class of piecewise functions with similar properties? The answer is the class of {\Chef{} splines}.\footnote{\Chefshort{} refers to the famous Russian mathematician and can be alternatively transliterated from the Russian writing of the name as Chebysheff, Chebyshev, Chebychov, Chebyshov, 
Tchebychev, Tschebyschev, Tschebyschef, Tschebyscheff, etc. We follow the traditional French transcription, in the footsteps of the monograph by \citet{Schumaker:2007}.} Besides algebraic polynomial splines, it contains exponential and trigonometric splines, and is closely related to null-spaces of linear differential operators \cite{Schumaker:2007}.

The term \emph{\Chef{} splines} was coined by \citet{KarlinZ:1966}. The authors considered functions belonging piecewise to a $(p+1)$-dimensional linear space $\ECT{\p}{}$, spanned by an extended complete \Chefshort{} system (ECT-system; see Section~\ref{sec:ECT-spaces-def}), and discussed their smoothness and approximation properties. 
We will refer to $p$ as the degree, in analogy with the polynomial splines. 
\Chef{} B-splines, among other properties, were established by \citet{Karlin:1968}.
It is impossible to give a complete account of the numerous articles related to this topic, so we do not want to make any attempt.
Noteworthy techniques for their construction and analysis are generalized divided differences \cite{Lyche:1985}, generalized polar forms \cite{Pottmann:1993}, generalized de Boor--Fix dual functionals \cite{Barry:1996} and repeated integration \cite{BisterP:1997}. 
\Chef{} B-splines possess all the characterizing properties of the classical polynomial B-splines. 
We refer the reader to the monograph by \citet{Schumaker:2007} and the survey article by \citet{Lyche:2019} for more details on \Chef{} splines and \Chef{} B-splines.

\Chef{} splines find applications in data approximation/interpolation \cite{KochL1993},
geometric modeling \cite{BeccariCM:2019,Mazure:2011a} and
signal processing \cite{Unser:2005,UnserB:2005}. Because of their relation to null-spaces of differential operators, they also offer a lot of opportunity in the context of isogeometric analysis, a spline paradigm for the numerical solution of differential problems \cite{Cottrell:2009}. Thanks to their structural similarities, \Chef{} B-splines are 
plug-to-plug compatible with classical polynomial B-splines, so they can be potentially easily incorporated in any software library supporting polynomial B-splines to enrich its capability.

A particularly interesting subclass of \Chef{} B-splines are the so-called \emph{generalized polynomial B-splines}, introduced by \citet{KvasovS:1999}. 
They can be seen as the minimal extension of (algebraic) polynomial B-splines towards the wide variety of \Chef{} B-splines, with a small selection of shape parameters. The fine-tuning of these parameters generally results in a gain from the accuracy point of view, compared with polynomial B-splines. In addition, suitable choices of such spaces --- including algebraic polynomial and exponential/trigonometric functions --- allow for an exact representation through (almost) arc-length parameterization of profiles of salient interest in applications, such as conic sections and helices. 
These are prominent features for geometric modeling \cite{Fang:2010,Wang:2008} and isogeometric analysis \cite{Aimi:2017,ManniPS:2011,ManniRS:2015,ManniRS:2017}.
A stable but costly method (based on convolution) to approximately evaluate generalized polynomial B-splines on uniform knots was proposed by
\citet{RomanMS:2017}.

It was shown by \citet{Nurnberger:1983,Nurnberger:1984} that many properties of ordinary \Chef{} splines carry over to certain generalized \Chef{} splines, in the sense that pieces can be drawn from different ECT-spaces of different dimensions. Under quite mild assumptions, such splines can be represented in terms of a B-spline-like basis. This basis is called \emph{generalized \Chef{} B-spline basis} or also \emph{multi-degree \Chef{} B-spline basis} to reflect better the analogy with polynomial splines. We will follow the latter terminology and refer to these basis functions as MDTB-splines. Their properties were studied more recently by \citet{Buchwald:2003} and \citet{Hiemstra:2020}.
The (algebraic) polynomial subclass of such multi-degree B-splines were explored in the context of geometric modeling \cite{Beccari:2017} and isogeometric analysis \cite{Toshniwal:2017polar}; these splines are called \emph{polynomial MDB-splines} or just \emph{MDB-splines}.

Unfortunately, despite their theoretical interest and applicative potential, MDTB-splines have not gained much attention in practice. The reason behind this is that MDTB-splines are generally difficult to compute. Classical approaches based on generalized divided differences, Hermite interpolation or repeated integration are computationally expensive and/or numerically unstable. An important step forward was recently made by \citet{Hiemstra:2020}; the authors proposed a construction based on a so-called \emph{multi-degree spline extraction operator} that represents MDTB-splines as linear combinations of local \Chef{} Bernstein functions. The local \Chef{} Bernstein functions form a basis of the local ECT-spaces involved in the definition of the MDTB-splines. In the polynomial case, these are nothing but the classical Bernstein polynomial basis functions.
The same type of extraction operator was already earlier investigated by \citet{Toshniwal:2017polar,Toshniwal:2020} and \citet{Speleers:2019} for dealing with the subclass of polynomial MDB-splines. A similar idea has also been pursued by \citet{Beccari:2021} for computing polynomial MDB-splines.

Here, we present an object-oriented \Matlab{} toolbox to construct and manipulate MDTB-splines whenever they exist. The key ingredient is the extraction operator discussed above. The toolbox is a continuation and extension of the \Matlab{} toolbox developed by \citet{Speleers:2019} for dealing with polynomial MDB-splines. The toolbox supports MDTB-splines whose pieces belong to ECT-spaces that are null-spaces of constant-coefficient linear differential operators. The computation of the corresponding \Chef{} Bernstein functions is inspired by the state-of-the-art implementation from the C++ library of \citet{Roth:2019} for the general ECT-space setting, but also relies on more efficient and more robust routines for certain specialized ECT-spaces (polynomial and generalized polynomial spaces of exponential and trigonometric type).

To the best of our knowledge, no general-purpose software library is nowadays available to work with \Chef{} splines, also considering the more restricted case where the local ECT-spaces have all the same degree or even where these local spaces are taken all the same. \Chef{} Bernstein functions and curves have been addressed by \citet{Roth:2019}, but spline curves are preferred in practice as they combine more local control of the shapes with built-in higher smoothness globally. Note that splines of lower degrees are usually employed for geometric modeling (so avoiding a source of ECT-space instabilities; see Section~\ref{sec:instability}). In this perspective, the \Matlab{} toolbox may unlock \Chef{} splines for a wide audience, and help pushing them from an elegant theoretical extension of polynomial splines towards a mainstream practical tool.

The remainder of the article is organized as follows. In Section~\ref{sec:ECT-spaces} we introduce the notion of ECT-space and show how to define \Chef{} Bernstein functions in such space. We detail in particular the important large class of ECT-spaces that are null-spaces of constant-coefficient linear differential operators. In Section~\ref{sec:spline-spaces} we focus on multi-degree \Chef{} spline spaces, and give a (theoretical) recursive definition of MDTB-splines. We also describe a knot insertion procedure that represents a set of MDTB-splines in terms of another set of MDTB-splines of lower smoothness. This procedure will form the foundation of the practical computation of MDTB-splines elaborated in Section~\ref{sec:computation}; it gives rise to an extraction operation that maps local \Chef{} Bernstein functions to the MDTB-spline basis of interest. In Section~\ref{sec:implementation} we discuss some practical implementation aspects and  review the general structure of the object-oriented \Matlab{} toolbox.
Section~\ref{sec:examples} illustrates the \Matlab{} toolbox with a selection of numerical examples, and we highlight certain pitfalls of working with ECT-spaces. We end in Section~\ref{sec:conclusion} with some concluding remarks.

\section{Extended Complete \Chefshort{} Spaces} \label{sec:ECT-spaces}

In this section, we define notation for ECT-spaces and recall some of their main properties. We also discuss an important basis for such spaces, the so-called \Chef{} Bernstein basis. We refer the reader to the survey works of \citet{Schumaker:2007} and \citet{Lyche:2019} for more details.

\subsection{ET-Spaces and ECT-Spaces} \label{sec:ECT-spaces-def}

We start by defining two important classes of \Chefshort{} spaces on a real interval $\interval$. 

\begin{definition}[Extended \Chefshort{} Space]
\label{def:ET}
Given an interval $\interval$, a space $\ECT{\p}{}(\interval)\subset C^{\p}(\interval)$ of dimension $\p+1$ ($\p\geq0$) is an {extended \Chefshort{} \mbox{(ET-)} space} on $\interval$ if any Hermite interpolation problem with $\p+1$ data on $\interval$ has a unique solution in $\ECT{\p}{}(\interval)$.
In other words, for any positive integer $m$,
let ${\bar x}_1,\ldots, {\bar x}_m$ be distinct points in $\interval$ and
let $d_1,\ldots,d_m$ be non-negative integers such that $\p+1=\sum_{i=1}^m (d_i+1)$. Then,
for any set $\{f_{i,j}\in\RR\}_{i=1,\ldots,m,\, j=0,\ldots,d_i}$
there exists a unique $g\in\ECT{\p}{}(\interval)$ such that
\begin{equation*}
D_{}^j g({\bar x}_i)=f_{i,j}, \quad i=1,\ldots,m, \quad j=0,\ldots,d_i.
\end{equation*}
\end{definition}

\begin{definition}[Extended Complete \Chefshort{} Space] \label{def:ECT}
Given an interval $\interval$, a space $\ECT{\p}{}(\interval)\subset C^{\p}(\interval)$ of dimension $p+1$ is an {extended complete \Chefshort{} (ECT-) space} if there exists a basis $\{\uu_0, \ldots, \uu_\p\}$ of $\ECT{\p}{}(\interval)$ such that every subspace $\langle \uu_0,\ldots,\uu_k\rangle$ is an ET-space on $\interval$ for $k=0,\ldots,\p$. The basis $\{\uu_0, \ldots, \uu_\p\}$ is called an {ECT-system}.
\end{definition}

A $(p+1)$-dimensional subspace of $C^{\p}(\interval)$ is an ECT-space on $\interval$ if and only if there exists a basis $\{\uu_0, \ldots, \uu_\p\}$ such that their \emph{Wronskian determinants} are positive:
\begin{equation*}
W[\uu_0,\ldots,\uu_k](x):=\det\begin{bmatrix}
 \uu_0(x) & \uu_1(x) & \cdots & \uu_k(x) \\
 D\uu_0(x) & D\uu_1(x) & \cdots & D\uu_k(x) \\
 \vdots & \vdots & & \vdots \\
 D^k\uu_0(x) & D^k\uu_1(x) & \cdots & D^k\uu_k(x) \\ 
\end{bmatrix}>0,\quad k=0,\ldots,\p,
\end{equation*}
 for all $x\in \interval$.
This basis forms an ECT-system and gives rise to a set of positive \emph{weight functions} defined by
\begin{equation} \label{eq:weights}
w_{j}(x) := \frac{W[\uu_0,\ldots,\uu_j](x)W[\uu_0,\ldots,\uu_{j-2}](x)}{\bigl(W[\uu_0,\ldots,\uu_{j-1}](x)\bigr)^2}, \quad j=0,\ldots,p,
\end{equation}
with the convention that $W[\emptyset]:=1$.
Conversely, any set of positive functions $w_{j} \in C^{\p-j}(\interval)$, $j=0,1,\ldots,\p$, generates the following ECT-system:
\begin{equation}
\label{eq:gen-powers}
\begin{cases}
	\uu_0(\x) := w_0(\x), \\
	\uu_1(\x) := w_0(\x) \int_{\aa}^\x w_1(\y_1) \dy_1, \\
	\hspace*{0.95cm} \vdots \\
	\uu_\p(\x) := w_0(\x) \int_{\aa}^\x w_1(\y_1) \int_{\aa}^{\y_1} \cdots   \int_{\aa}^{y_{\p-1}} w_\p(\y_\p) \dy_{\p} \cdots \dy_1,
\end{cases}
\end{equation}
for any fixed point $\aa\in\interval$. The functions $\uu_0, \ldots, \uu_\p$ in \eqref{eq:gen-powers} are called \emph{generalized powers}. 
From a practical point of view, it is often desired that the space $\ECT{\p}{}(\interval)$ contains constants. This is achieved with the choice $w_0=1$. 

\begin{example} \label{ex:poly}
The space $\PP_\p:=\langle 1, x, \ldots,x^\p\rangle$ of algebraic polynomials is an ECT-space on any interval of the real line.
It can be regarded as the span of the ECT-system
\begin{equation} \label{eq:taylor-basis}
  \left\{1,{x-\aa},\frac{(x-\aa)^2}{2},\ldots,\frac{(x-\aa)^\p}{\p!}\right\},
\end{equation}
for any fixed point $\aa\in\RR$. Indeed, the Wronskian determinants of this system are all equal to one. The functions in \eqref{eq:taylor-basis} form the classical Taylor basis for algebraic polynomials. They can be generated by the weight functions $w_0=\cdots=w_\p=1$ according to \eqref{eq:gen-powers}.
\end{example}

\begin{remark}\label{rmk:weights_factor}
A given ECT-space can be identified by different sets of weight functions; see \citet{Lyche:2019} for details and examples. In particular, it is easy to see that
the two weight systems
\begin{equation*}
w_0,\ldots,w_\p \quad \text{and} \quad K_0w_0,\ldots,K_\p w_\p,
\end{equation*}
where $K_0,\ldots, K_\p$ are positive constants, identify the same ECT-space. A constructive procedure for finding all weight systems associated with a given ECT-space on a bounded closed interval is described by \citet{Mazure:2011b}.
\end{remark}

\begin{remark}
From Definition~\ref{def:ECT} it is clear that an ECT-space of dimension $p+1$ on $J$ is an ET-space of dimension $p+1$ on $J$. The converse is not true in general. However, if $\interval$ is a bounded closed interval, then any ET-space of dimension $\p+1$ on $\interval$ is an ECT-space on $\interval$; see \citet{Mazure:2007}. In the context of \Chef{} spline spaces, the scope of this article, we are only interested in bounded closed intervals, so both notions are interchangeable. Further on, we use the notion ECT even if it can be weakened to ET.
\end{remark}

\subsection{\Chef{} Bernstein Functions}
We now set $\interval:=[\x_0,\x_1]$ with $\x_0<\x_1$. Instead of working with the generalized power basis, an alternative basis is formed by the so-called \emph{\Chef{} Bernstein functions} associated with the ECT-space $\ECT{\p}{}(\interval)$. They are denoted with $\bs_{j,\p}$, $j=0, \ldots, \p$ and can be defined recursively as follows. Let $w_j$, $j=0, \ldots, \p$ be positive weight functions generating $\ECT{\p}{}(\interval)$ 
and we assume $w_0=1$. For $\q=0,\ldots, \p$ and $j=0,\ldots,\q$, the function $\bs_{j,\q}$ is defined at $\x \in [\x_0, \x_1]$ as
\begin{equation} \label{eq:rec-bern-0}
\bs_{0,0}(\x) := w_{\p}(\x),
\end{equation}
and
\begin{equation}  \label{eq:rec-bern-q}
\bs_{j,\q}(\x) := w_{\p-\q}(\x)\cdot\begin{cases}
\displaystyle 1-\int_{\x_{0}}^\x \dfrac{\bs_{0,\q-1}(\y)}{\bbint_{0,\q-1}}\dy, & j=0,
\\[0.35cm]
\displaystyle \int_{\x_{0}}^\x \biggl[\dfrac{\bs_{j-1,\q-1}(\y)}{\bbint_{j-1,\q-1}} - \frac{\bs_{j,\q-1}(\y)}{\bbint_{j,\q-1}}\biggr]\dy, &  0< j < q,
\\[0.35cm]
\displaystyle \int_{\x_{0}}^\x \dfrac{\bs_{\q-1,\q-1}(\y)}{\bbint_{\q-1,\q-1}}\dy, & j=\q,
\end{cases}\quad \q>0,
\end{equation}
where
\begin{equation*}
\bbint_{j,\q-1} := \int_{\x_{0}}^{\x_1} \bs_{j, \q-1}(\y)\dd\y.
\end{equation*}
The Bernstein functions $\bs_{0,\p},\ldots,\bs_{\p,\p}$ are non-negative, form a partition of unity, and enjoy the following end-point conditions:
\begin{equation} \label{eq:endpoint}
\begin{aligned}
\bs_{0,\p}(\x_{0}) &=1, \quad
  D^{k}\bs_{j,\p}(\x_{0})=0, \quad k=0,\ldots,j-1, \\
\bs_{\p,\p}(\x_{1})&=1, \quad
  D^{k}\bs_{j,\p}(\x_{1})=0, \quad k=0,\ldots,\p-j-1.
\end{aligned}
\end{equation}
Moreover, they are a basis of the space $\ECT{\p}{}(\interval)$.

\begin{example}\label{ex:poly-bernstein}
When dealing with algebraic polynomials, see Example~\ref{ex:poly}, the \Chef{} Bernstein functions are nothing but the classical Bernstein polynomials, which can be explicitly expressed as
\begin{equation*}  
\bs_{j,\p}(\x)=\frac{\p!}{j!(\p-j)!}\left(\frac{\x-\x_0}{\x_1-\x_0}\right)^{j}\left(\frac{\x_1-\x}{\x_1-\x_0}\right)^{\p-j},
\quad j=0, \ldots, \p.
\end{equation*}
\end{example}

\begin{remark}\label{rmk:existence-Bernstein}
The definition of the \Chef{} Bernstein basis in a given ECT-space $\ECT{\p}{}(\interval)$ requires the existence of a set of positive weight functions $w_j$, $j=0, \ldots, \p$, generating $\ECT{\p}{}(\interval)$ such that $w_0=1$. This is guaranteed if and only if the derivative space of $\ECT{\p}{}(\interval)$ is an ECT-space. 
\end{remark}

\begin{remark}\label{rmk:Hermite-Bernstein}
Instead of using the recurrence relation \eqref{eq:rec-bern-0}--\eqref{eq:rec-bern-q}, each Bernstein function $\bs_{j,\p}$ can also be computed by solving the following Hermite interpolation problem in the space $\ECT{\p}{}(\interval)$: for $j=0$,
\begin{equation*}
\bs_{0,\p}(\x_{0})=1, \quad
  D^{k}\bs_{0,\p}(\x_{1})=0, \quad k=0,\ldots,\p-1,
\end{equation*}
and for $j=1,\ldots, \p$,
\begin{equation} \label{eq:Hermite-Bernstein}
\begin{aligned}
D^{k}\bs_{j,\p}(\x_{0})&=0, \quad  k=0,\ldots,j-1, \quad
D^{k}\bs_{j,\p}(\x_{1})=0, \quad k=0,\ldots,\p-j-1, \\
D^{j}\bs_{j,\p}(\x_{0})&=-\sum_{k=0}^{j-1}D^{j}\bs_{k,\p}(\x_{0}).
\end{aligned}
\end{equation}
Since $\ECT{\p}{}(\interval)$ is an ECT-space, this interpolation problem has a unique solution; see Definition~\ref{def:ET}. Note that the conditions in \eqref{eq:Hermite-Bernstein} require that the Bernstein functions $\bs_{j,\p}$ are computed sequentially from $j=0$ to $j=p$. The order can be reversed by employing the alternative conditions
\begin{equation*}
  \bs_{p,\p}(\x_{1})=1, \quad 
  D^{\p-j}\bs_{j,\p}(\x_{1}) =-\sum_{k=j+1}^{\p}D^{\p-j}\bs_{k,\p}(\x_{1}),
  \quad j=0,\ldots,\p-1.
\end{equation*}
Any convenient basis of $\ECT{\p}{}(\interval)$ can be used to represent the Bernstein functions.
\end{remark}

\subsection{A Large Class of ECT-Spaces} \label{sec:ECT-Np}

Let $\cL_p$ be the {linear differential operator} defined by
\begin{equation}\label{eq:Lp}
\cL_\p f := D^{\p+1}f+\sum_{j=0}^{\p} a_j D^jf,  \quad f\in C^{\p+1}(\interval),
\end{equation}
with constant coefficients $a_j\in\RR$ and $\interval:=[\x_0,\x_1]$.
Any operator of the form \eqref{eq:Lp} is uniquely identified by its null-space $\NN_\p$. A fundamental set of solutions, forming a basis of $\NN_\p$, can be generated through the (higher-order) roots of the characteristic polynomial
\begin{equation}\label{eq:Lp-pol}
  \croot^{\p+1}+\sum_{j=0}^{\p} a_j \croot^j, \quad \croot\in\CC,
\end{equation}
associated with the differential operator in \eqref{eq:Lp}. Let $\croot=\alpha+\ii \beta$ be a root of order $\mu$ ($\mu\geq1$) of the polynomial in \eqref{eq:Lp-pol} for some $\alpha,\beta\in\RR$ and $\ii:=\sqrt{-1}$. Then, this root generates the following fundamental subspace:
\begin{itemize}
\item if $\beta=0$, then 
\begin{equation*}
  \bigl\langle \x^i \ee^{\alpha\x}:i=0,\ldots,\mu-1\bigr\rangle \subseteq \NN_\p;
\end{equation*}
\item if $\beta\neq0$, then the complex conjugate of $\croot$ is also a root of order $\mu$, and 
\begin{equation*}
  \bigl\langle \x^i \ee^{\alpha\x}\cos(\beta\x),\x^i \ee^{\alpha\x}\sin(\beta\x):i=0,\ldots,\mu-1\bigr\rangle \subseteq \NN_\p.
\end{equation*}
\end{itemize}
The fundamental subspaces related to different (non-conjugate) roots are disjoint, and all together they span the full null-space.
Note that $\NN_\p$ is translation-invariant. In order to ensure that constants belong to $\NN_\p$, we have to assume that $\croot=0$ is at least a first-order root of the characteristic polynomial \eqref{eq:Lp-pol}.

The null-space $\NN_\p$ is an ECT-space on $\interval$ if and only if there exist positive weight functions $w_j\in C^{\p-j+1}(\interval)$, $j=0,\ldots,\p$ such that
\begin{equation}\label{eq:Lp-factor}
 \cL_\p f = w_0\cdots w_\p D_p \cdots D_0f,
\end{equation}
where
\begin{equation*}
D_j f := D\left(\frac{f}{w_j}\right), \quad j=0,\ldots,\p;
\end{equation*}
see \citet[Chapter~3]{Coppel:1971}. Furthermore, $\NN_\p$ is always an ECT-space on intervals of sufficiently small length. 
The so-called \emph{critical length}, $\ell_\p>0$, is the supremum of the range of lengths of the intervals on which the space is ECT. 
The critical length can be bounded from below as
\begin{equation*}
\ell_\p\geq \pi/M\in(0,+\infty],
\end{equation*}
where $M\geq0$ is the maximum of the imaginary parts of all roots of the characteristic polynomial. In view of Remark~\ref{rmk:existence-Bernstein}, the existence of a \Chef{} Bernstein basis in $\NN_\p$ is guaranteed if and only if the derivative space of $\NN_\p$ is ECT, and hence the critical length of the derivative space is also of practical importance. This length is denoted with $\ell'_\p$ and often called \emph{critical length for design}.
We refer the reader to \citet{Carnicer:2003} for a theoretical study on critical lengths of such ECT-spaces and to \citet{BeccariCM:2020} for a numerical study.

\begin{example} \label{ex:nullspace-poly}
From Example~\ref{ex:poly} we know that the polynomial space $\PP_\p$ is an ECT-space generated by the weight functions $w_0=\cdots=w_\p=1$. This space is the null-space of the linear differential operator $\cL_\p f=D^{\p+1}f$, in agreement with \eqref{eq:Lp-factor}, and its critical length is $\ell_\p=+\infty$.
\end{example}

\begin{example} \label{ex:nullspace-exp}
Consider the linear differential operator
\begin{equation*}
\cL_\p f = D^{\mu_0}(D-\alpha_1)^{\mu_1}\cdots(D-\alpha_m)^{\mu_m}f, 
\end{equation*}
with distinct real values $\alpha_k\neq0$, $k=1,\ldots,m$, and positive integer values $\mu_k$, $k=0,\ldots,m$ such that $\sum_{k=0}^m \mu_k=\p+1$. Its null-space can be described by means of the fundamental set of solutions as
\begin{equation}\label{eq:nullspace-exp}
\NN_\p = \bigl\langle 1,x,\ldots,x^{\mu_0-1},
\ee^{\alpha_1 \x},x\ee^{\alpha_1 \x},\ldots,\x^{\mu_1-1} \ee^{\alpha_1 \x},
\ldots,
\ee^{\alpha_m \x},\x\ee^{\alpha_m \x},\ldots,\x^{\mu_m-1} \ee^{\alpha_m \x}\bigr\rangle.  
\end{equation}
It is easily verified that
\begin{equation*}
W[\ee^{\alpha \x},x\ee^{\alpha \x},\ldots,\x^{\mu-1} \ee^{\alpha \x}](x)
= \ee^{\alpha \mu\x} \prod_{k=0}^{\mu-1} k!>0,
\end{equation*}
for any $\alpha,\x\in\RR$ and $\mu\in\NN$. More generally, the set of functions in \eqref{eq:nullspace-exp} forms an ECT-system on any interval in $\RR$, and can be used to compute a set of positive weight functions according to \eqref{eq:weights}. This implies that the space in \eqref{eq:nullspace-exp} is an ECT-space and its critical length is $\ell_\p=+\infty$. The derivative space of $\NN_\p$ in \eqref{eq:nullspace-exp} is of the same form as $\NN_\p$, so its critical length is $\ell'_\p=+\infty$.
\end{example}

\begin{example} \label{ex:nullspace-trig}
Consider the linear differential operator
\begin{equation*}
\cL_\p f = D^{\p-1}(D^2+\beta^2)f, 
\end{equation*}
for some real value $\beta\neq0$ and $\p\geq1$. Its null-space can be described by means of the fundamental set of solutions as
\begin{equation}\label{eq:nullspace-trig}
\NN_1 = \bigl\langle \cos(\beta\x),\sin(\beta\x)\bigr\rangle,  
\quad
\NN_\p = \bigl\langle 1,x,\ldots,x^{\p-2},\cos(\beta\x),\sin(\beta\x)\bigr\rangle,
\quad \p\geq2.
\end{equation}
Such space is sometimes called {cycloidal space}. It is an ECT-space on finite intervals of sufficiently small length. Assuming $\beta=1$, it has been shown by \citet{Carnicer:2003} that $\ell_\p \leq 2\pi\lfloor \p/2 \rfloor$ for $\p\geq2$ and in particular that $\ell_1=\pi$, $\ell_2=\ell_3=2\pi$, and $\ell_4=\ell_5\approx8.9868$. More generally, the critical lengths are connected to roots of Bessel functions \cite{CarnicerMP:2017}.
The derivative space of $\NN_\p$ in \eqref{eq:nullspace-trig} is of the same form as $\NN_\p$ and its critical length is found to be $\ell'_\p=\ell_{\p-1}$ for $\p\geq2$. 
\end{example}

\section{Multi-Degree \Chef{} Spline Spaces} \label{sec:spline-spaces}

We are interested in spaces of piecewise functions, whose pieces belong to ECT-spaces and are glued together in a certain smooth way. We show how a B-spline-like basis can be defined for such spaces, the so-called MDTB-spline basis, with similar characteristics to the classical polynomial B-spline basis. We summarize the main properties and follow the notation from \citet{Hiemstra:2020}; see also \citet{Buchwald:2003} and \citet{Nurnberger:1984}.

\subsection{Multi-Degree \Chef{} Splines}

Pieces of our splines shall be drawn from arbitrary ECT-spaces of possibly different dimensions. Consider a partitioning of the interval $[\a,\b] \subset \RR$ into a sequence of break points,
\begin{equation*}
  \domain := \{\a =: \x_0 < \x_1 <  \cdots < \x_{\nelms-1} < \x_{\nelms} := \b\}.
\end{equation*}
Furthermore, we set $\interval_i := [\x_{i-1}, \x_{i})$, $i=1, \ldots, \nelms-1$, and $\interval_\nelms := [\x_{\nelms-1}, \x_{\nelms}]$. We also define an ECT-space of dimension $\p_i+1$ on each closed interval $[\x_{i-1}, \x_{i}]$, $i=1, \ldots, \nelms$:
\begin{equation*}
  \ECT{\p_i}{(i)} := \bigl\langle\uu^{(i)}_0, \ldots , \uu^{(i)}_{\p_i}\bigl\rangle, \quad \uu^{(i)}_j \in C^{\p_i}([\x_{i-1}, \x_{i}]), \quad j=0, \ldots, \p_i,
\end{equation*}
where $\uu^{(i)}_0,\ldots,\uu^{(i)}_{\p_i}$ are generalized powers defined in terms of positive weight functions $w_j^{(i)}\in C^{\p_i-j}([\x_{i-1},\x_{i}])$, $j=0,\ldots,\p_i$ as in \eqref{eq:gen-powers}.
%
Collectively, these local ECT-spaces span the following global space:
\begin{equation}\label{eq:spline-space-discont}
\splSpacep := \bigl\{ { s: [\a, \b] \rightarrow \RR :}
  \left. s \right|_{\interval_i} \in \ECT{\p_i}{(i)},\; i=1,\ldots, \nelms \bigr\}.
\end{equation}
Then, by prescribing the smoothness at the break points we can define the corresponding space of \Chef{} splines as follows.

\begin{definition}[MDT-Spline Space] \label{def:spline-space}
Given the sets of integers $\bp:=\{\p_1,\ldots,\p_\nelms\}$ and
\begin{equation} \label{eq:smoothness}
\br := \bigl\{\r_i \in \ZZ  : -1 \leq \r_i \leq \min\{\p_{i}, \p_{i+1}\}, \; i=1,\ldots , \nelms-1, \; \r_0 = \r_{\nelms} = -1 \bigr\},
\end{equation}
 we define
\begin{equation*} 
\splSpacerp := \bigl\{  s \in \splSpacep :  D^{j}_{-}s(\x_i) = D^{j}_{+}s(\x_i), \; j = 0,\ldots, \r_i \text{ and } i = 1, \ldots, \nelms-1 \bigr\}.
\end{equation*}
This space is called multi-degree \Chef{} spline  (MDT-spline) space.
\end{definition}

The value $\r_i$ represents the smoothness at break point $\x_i$, $i=1,\ldots,m-1$. All smoothness conditions are linearly independent because the functions $\{\uu^{(i)}_0,\ldots,\uu^{(i)}_{\p_i}\}$ on each interval $\interval_i$ form an ECT-system. 
Hence, the dimension of $\splSpacerp$ is given by
\begin{equation}\label{eq:spline-dimension}
\N := \sum_{i=0}^{\nelms-1} {(\p_{i+1}-\r_i)}=\sum_{i=1}^{\nelms}(\p_i-\r_i).
\end{equation}

When considering different ECT-spaces on different intervals, the construction of spline spaces equipped with the same properties as classical polynomial splines (including a B-spline-like basis) requires constraints on the various ECT-spaces.

\begin{definition}[Admissible Weights]
\label{def:weight-assumption}
The weight systems $\{ w_j^{(i)}: j=0,\ldots,\p_i \}$ generating the ECT-spaces $\ECT{\p_i}{(i)}$, $i=1,\ldots,m$, are {admissible} for the space $\splSpacerp$ if 
\begin{equation*}
D_{-}^lw_j^{(i)}(x_i)=D_{+}^lw_j^{(i+1)}(x_i), \quad l=0,\ldots, \r_i-j,
\end{equation*}
for $i=1,\ldots,\nelms-1$ and $j=0,\ldots,\r_i$.
Moreover, it is assumed that $w_0^{(i)}=1$ for $i=1,\ldots,\nelms$.
\end{definition}

\begin{remark}\label{rmk:weight-assumption}
Dealing with admissible weights gives only a sufficient condition for obtaining \Chef{} splines equipped with a B-spline-like basis; see \citet{Buchwald:2003}. The simplicity of this condition and the fact that it embraces relevant classes of \Chef{} splines motivate our choice. 
We refer the reader to \citet{Mazure:2018} for explicit necessary and sufficient conditions for smoothly gluing together ECT-spaces of dimension $5$.
\end{remark}

\begin{remark}\label{rmk:periodicity}
  Definition~\ref{def:spline-space} can be easily extended to incorporate periodicity. In this case, periodic continuity constraints need to be imposed and we set $\r_0 = \r_{\nelms} = \r_{{\rm per}}$ for some value of $\r_{{\rm per}}\geq0$. For simplicity of exposition, we do not consider this extension in the following.
\end{remark}

\subsection{Multi-Degree \Chef{} B-Splines}

We now introduce basis functions for the MDT-spline space $\splSpacerp$ that possess all the characterizing properties of classical polynomial B-splines. We call the corresponding functions \emph{multi-degree \Chef{} B-splines (MDTB-splines)} to stress the fact that ECT-spaces of different dimensions can be employed on different intervals, in analogy with the polynomial MDB-splines considered in \citet{Speleers:2019} and \citet{Toshniwal:2020}.

The construction and analysis of MDTB-splines can be eased by considering two knot vectors,
\begin{align}
\mbf{\lknot}{} &:= (\lknot_k )_{k=1}^{\N}:= (\;
    \underbrace{\x_0,\; \ldots,\;  \x_0}_{\p_{1} - \r_0 \text{ times}},\;  \ldots,\;
    \underbrace{\x_{i},\;  \ldots,\;  \x_{i}}_{\p_{i+1} - \r_i \text{ times}},\;   \ldots,\;
    \underbrace{\x_{\nelms-1},\;  \ldots,\;  \x_{\nelms-1}}_{\p_{\nelms} - \r_{\nelms-1} \text{ times}}
    \; ), \label{eq:knots-u}\\
\mbf{\rknot}{} &:= (\rknot_k )_{k=1}^{\N}:= (\;
    \underbrace{\x_1,\;  \ldots,\;  \x_1}_{\p_1 - \r_1 \text{ times}},\; \ldots,\;
    \underbrace{\x_{i},\;  \ldots,\;  \x_{i}}_{\p_i - \r_i \text{ times}},\; \ldots ,\;
    \underbrace{\x_\nelms,\;  \ldots,\;  \x_\nelms}_{\p_{\nelms} - \r_{\nelms} \text{ times}}
    \; ). \label{eq:knots-v}
\end{align}

Assume there exist admissible weights for the space $\splSpacerp$. The set of MDTB-splines $\{\BS_k: k=1,\ldots,\N\}$ can be computed through an integral recurrence relation that is very similar to the one of the \Chef{} Bernstein functions in \eqref{eq:rec-bern-0}--\eqref{eq:rec-bern-q}. To this end, we set $\p := \max_{1\leq i \leq\nelms} \p_i$, and we define a global set of weight functions $\{w_j: j=0,\ldots,\p\}$ by
\begin{equation*} 
w_j(\x) :=
\begin{cases}
    w_{j}^{(i)}(\x), & j\leq \p_i, \\
    0, & \text{otherwise},
\end{cases}
\quad \x \in [\x_{i-1}, \x_{i}), \quad i=1,\ldots,\nelms.
\end{equation*}
Then, the MDTB-splines $\BS_k = \BS_{k,\p}$, $k=1, \ldots, \N$, can be defined recursively as follows. 
For $\q=0,\ldots, \p$ and $k=\p-\q+1,\ldots,\N$, the spline $\BS_{k,\q}$ is supported on the interval $[\lknot_k,\rknot_{k - \p+\q}]$, and is defined at $\x \in [\a,\b)$ as 
\begin{equation} \label{eq:rec-global-0}
\BS_{k, 0}(\x) :=
\begin{cases}
    w_{\p}(\x), & \x \in [\x_{i-1}, \x_{i}),\\
    0, & \text{otherwise},
\end{cases}
\end{equation}
and
\begin{equation}\label{eq:rec-global-q}
\BS_{k,\q}(\x) :=
    w_{\p-\q}(\x)\int_{\a}^\x \biggl[\dfrac{\BS_{k,\q-1}(\y)}{\bspint_{k,\q-1}} - \dfrac{\BS_{k+1,\q-1}(\y)}{\bspint_{k+1,\q-1}}\biggr] \dd \y, \quad \q > 0,
\end{equation}
where
\begin{equation*} 
\bspint_{j,\q-1} := \int_{\a}^{\b} \BS_{j, \q-1}(\y)\dd\y.
\end{equation*}
In the above we assumed that any undefined $\BS_{j, \q-1}$ with $j<\p-\q+2$ or $j>\N$ must be regarded as the zero function,
and we used the convention that if $\bspint_{j, \q-1}=0$ then
\begin{equation*} 
\int_{\a}^\x \dfrac{\BS_{j,\q-1}(\y)}{\bspint_{j,\q-1}}\dd\y :=
\begin{cases}
1, & \x \geq \lknot_{j} \text{~and~} j\leq \N,\\
0, & \text{otherwise}.
\end{cases}
\end{equation*}
At the right end point $\b$, the spline  $\BS_{k,\q}$ is defined by taking the limit from the left, that is $\BS_{k,\q}(\b):=\lim_{\x\rightarrow \b,\x<\b}\BS_{k,\q}(\x)$. We refer the reader to \citet{Hiemstra:2020} for alternative definitions.

\begin{example} \label{ex:MDTB-Bernstein}
Consider the discontinuous MDT-spline space $\splSpacep$ in \eqref{eq:spline-space-discont}. Each of the corresponding MDTB-splines is supported on a single interval $[\x_{i-1},\x_{i}]$ for some $i$. Moreover, the non-zero MDTB-splines on $[\x_{i-1},\x_{i})$ coincide on this interval with the \Chef{} Bernstein functions associated with the ECT-space $\ECT{\p_i}{(i)}$.  
  Hence, the similarity between the definitions in \eqref{eq:rec-bern-0}--\eqref{eq:rec-bern-q} and \eqref{eq:rec-global-0}--\eqref{eq:rec-global-q} is not a coincidence.
\end{example}

The MDTB-spline basis enjoys several nice properties.

\begin{proposition} \label{pro:properties}
Assume there exist admissible weights for the space $\splSpacerp$. Then, the set $\{\BS_k: k=1, \ldots, \N\}$ is a basis of the space $\splSpacerp$, with the following properties:
\begin{itemize}
  \item local support:
  \begin{equation*}
  \supp(\BS_k)=[\lknot_k, \rknot_{k}], \quad k=1,\ldots,\N;
  \end{equation*}
  \item non-negative partition of unity:
  \begin{equation*}
    \BS_k(\x)\geq0, \quad k=1,\ldots,\N, \quad \sum_{k=1}^{\N}\BS_k(\x)=1, \quad x\in[\a,\b]; 
  \end{equation*}
  \item interpolation at the end points:
  \begin{align*} 
   &\BS_{1}(\a)=1, \quad \BS_{k}(\a)=0, \quad k=2,\ldots,\N, \\
   &\BS_{\N}(\b)=1, \quad \BS_{k}(\b)=0, \quad k=1,\ldots,\N-1.
  \end{align*}
\end{itemize}
\end{proposition}

These properties are of interest in both geometric modeling and isogeometric analysis; they make the set of MDTB-splines $\{\BS_k: k=1, \ldots, \N\}$ a very appealing basis for the space $\splSpacerp$ in those applications. However, using the recurrence relation in \eqref{eq:rec-global-0}--\eqref{eq:rec-global-q} for their construction is a computational nightmare! Not only one has to cope with a repeated calculation of integrals, but also one has to find a proper set of weight functions that is admissible for the MDT-spline space. 
Given the complexity of finding such weight functions in general (see, e.g., \citet{Lyche:2019} and references therein), it is wishful to avoid them in the computation at all. 
An alternative way of constructing MDTB-splines is based on knot insertion, where a new set of basis functions is computed from another set of basis functions. The idea is outlined in the following section.

\subsection{Knot Insertion}
Let us first observe that MDTB-splines possess super-smoothness at the knots, that is higher smoothness than the space $\splSpacerp$ requires. More precisely, if $\r_i<\min\{\p_i,\p_{i+1}\}$, then
there are only $\r_i+3$ successive MDTB-splines $\BS_k$ that have a jump in their $(\r_i+1)$-th order derivative at $\x = \x_i$, namely for $k = \lintsum(i), \ldots , \rintsum(i)+1$,
where
\begin{equation*}
\rintsum(i) := \sum_{j=0}^{i-1} (\p_{j+1} - \r_{j}), \quad
\lintsum(i) := \sum_{j=1}^{i} (\p_{j} - \r_j), \quad i=1,\ldots,\nelms-1.
\end{equation*}

Suppose now that the knot vectors $\mbf{\lknot}{} $ and $\mbf{\rknot}{}$ defined in \eqref{eq:knots-u}--\eqref{eq:knots-v} are obtained from other knot vectors $\tilde{\mbf{\lknot}{}} $ and $\tilde{\mbf{\rknot}{}} $ by inserting a single knot $\lknot = \rknot = \x_i \in (a, b)$, respectively, for some $0<i<\nelms$.
The related smoothness vector is easily deduced to be $\tilde{\br} = (\r_0, \ldots, \r_{i-1}, \r_{i}+1, \r_{i+1}, \ldots, \r_{\nelms})$; it is assumed to satisfy the same restrictions as in \eqref{eq:smoothness}, and so $\r_i+1\leq\min\{\p_i,\p_{i+1}\}$. Consequently, the spline space $\splSpacerptilde$ is a subspace of $\splSpacerp$ with one additional continuous derivative at $\x = \x_i$. Let $\{\tilde{\BS}_k:k=1,\ldots \N-1\}$ be the set of MDTB-splines of $\splSpacerptilde$.
Then,
\begin{equation}\label{eq:knotremoval}
\tilde{\BS}_k(\x) = \ca_k \BS_k(\x) + \cb_{k+1} \BS_{k+1}(\x),
\quad k=1,\ldots \N-1,
\end{equation}
where
\begin{itemize}
  \item[(i)] $\ca_{k} = 1$ and $\cb_{k+1}=0$ if $1\leq k < \lintsum(i)$;
  \item[(ii)] $\ca_{k} > 0$ and $\cb_{k+1} = - \ca_{k} \dfrac{l_{k}}{l_{k+1}} > 0$ if $\lintsum(i) \leq k \leq \rintsum(i)$;
  \item[(iii)] $\ca_{k} = 0$ and $\cb_{k+1}=1$ if $\rintsum(i) < k < \N$,
\end{itemize}
and
\begin{equation*}
l_k := {D^{\r_i+1}_{-}\BS_{k}(\x_i)-D^{\r_i+1}_{+}\BS_{k}(\x_i)}. \end{equation*}
Moreover, $\ca_{k}+\cb_{k}=1$ for $k=1,\ldots,\N$. This property implies that the coefficients in item (ii) can be computed in succession as follows:
\begin{equation} \label{eq:knotremoval-pattern} 
\begin{alignedat}{4} 
\ca_{\lintsum(i)}  = 1
\quad \rightarrow \quad
\cb_{\lintsum(i)+1} &= - \ca_{\lintsum(i)}  \frac{l_{\lintsum(i)}}{l_{\lintsum(i)+1}} \\[-0.1cm]
& \, \downarrow \\[-0.1cm]
\ca_{\lintsum(i)+1} &= 1 - \cb_{\lintsum(i)+1} 
\quad \rightarrow \quad&
\cb_{\lintsum(i)+2} &= - \ca_{\lintsum(i)+1} \frac{l_{\lintsum(i)+1}}{l_{\lintsum(i)+2}} \\[-0.1cm]
&&& \, \downarrow \\[-0.1cm]
&&& \quad\ddots\quad \\[-0.1cm]
&&& \quad\quad\quad\rightarrow\quad
 \cb_{\rintsum(i)+1} =1.
\end{alignedat}
\end{equation}
The relation in \eqref{eq:knotremoval} allows us to write a set of MDTB-splines in terms of another set of MDTB-splines of lower smoothness. Hence, low-smooth spaces can be used as a step-up to deal with high-smooth spaces. In this perspective, the discontinuous space $\splSpacep$ in \eqref{eq:spline-space-discont} is useful as starting point since the corresponding MDTB-splines can be locally computed as \Chef{} Bernstein functions; see Example~\ref{ex:MDTB-Bernstein}.

\section{Computational Aspects} \label{sec:computation}

As already discussed before, the computation of MDTB-splines through the integral recurrence relation in \eqref{eq:rec-global-0}--\eqref{eq:rec-global-q} is numerical challenging. Therefore, in this section, we describe an alternative, practical construction based on knot insertion \cite{Hiemstra:2020}. The construction relies on an extraction operator that represents all MDTB-splines as linear combinations of local \Chef{} Bernstein functions. For its practical implementation, we can closely follow the algorithmic procedure by \citet{Speleers:2019} developed for the specific case of polynomial MDB-splines.

\subsection{Computation of MDTB-Splines} \label{sec:computation-spline}

On the $i$-th interval $\interval_i$, we have $\p_i+1$ Bernstein functions $\bs_{0,\p_i}^{(i)},\ldots,\bs_{\p_i,\p_i}^{(i)}$ that span the local spline space $\ECT{\p_i}{(i)}$. 
In the first step, we extend them on the entire interval $[\a,\b]$ by defining them to be zero outside $\interval_i$.
More precisely, setting 
\begin{equation*}
  \nsum(0):=0, \quad 
  \nsum(i):=\sum_{j=1}^{i}(\p_{j}+1),\quad i=1,\ldots,\nelms, 
\end{equation*}
we define for $i\in\{1,\ldots,\nelms\}$, $j\in\{0,\ldots,\p_i\}$,
\begin{equation*}
  \bs_{\nsum(i-1) + j + 1}(\x) := 
  \begin{cases}
  \bs_{j,\p_i}^{(i)}(\x),& 
  \text{if } \x\in \interval_i,\\
  0, & \text{otherwise}.
  \end{cases}
\end{equation*}
For the sake of simplicity, we dropped the reference to the (local) degree in the notation.
From the properties of B-splines, it is clear that the functions $\bs_{1},\ldots,\bs_{\nsum(\nelms)}$ are linearly independent,
form a non-negative partition of unity, and span the space $\splSpacep$.
We arrange these basis functions in a column vector $\mbf{\bs}$ of length $\nsum(\nelms)$.

Now, we are looking for the set of MDTB-spline basis functions $\{\BS_1,\ldots,\BS_\N\}$ that span the smoother space $\splSpacerp$.
We arrange these basis functions in a column vector $\mbf{\BS}$ of length $\N$.
Since $\splSpacerp\subseteq\splSpacep$, we aim to construct a matrix $\mbf{H}$ of size $\N \times \nsum(\nelms)$ such that 
\begin{equation}\label{eq:multiDegBasisExt}
  \mbf{\BS} = \mbf{H}\,\mbf{\bs}.  
\end{equation}
To this end, we build continuity constraints at all break points corresponding to $\br$ and construct $\mbf{H}$ as their (left) null-space. For the computation of $\mbf{H}$, we can apply exactly the same algorithm as described by \citet[Section~3.2]{Speleers:2019} in the specific case of polynomial splines, thanks to the structural similarity between MDTB-splines and MDB-splines \cite{Hiemstra:2020}. For the sake of completeness and comprehension in our \Chef{} spline setting, we revisit the algorithm in the following.

Consider the $i$-th break point $\x_i$ for some $i\in\{1,\ldots,\nelms-1\}$.
Let $\mbf{K}^{(i,1)}$ be a matrix of size $(\p_i+1)\times(\r_i+1)$,
whose $\deriv$-th column is given by
\begin{equation} \label{eq:K_L}
    \begin{bmatrix}
    0 & \cdots & 0 & D_-^{\deriv-1}\bs_{\nsum(i)-\deriv+1}(\x_i) & \cdots & D_-^{\deriv-1}\bs_{\nsum(i)}(\x_i)
    \end{bmatrix}^T,
\end{equation}
and let $\mbf{K}^{(i,2)}$ be a matrix of size $(\p_{i+1}+1)\times(\r_i+1)$, whose $\deriv$-th column is given by
\begin{equation} \label{eq:K_R}
    \begin{bmatrix}
    -D_+^{\deriv-1}\bs_{\nsum(i)+1}(\x_i) & \cdots & -D_+^{\deriv-1}\bs_{\nsum(i)+\deriv}(\x_i) & 0 & \cdots & 0
    \end{bmatrix}^T.
\end{equation}
Note that the derivatives of the basis functions in the above matrices can be computed by evaluating the derivatives of the corresponding local Bernstein functions at the end points of their basic interval. This explains the triangular structure of both matrices; see \eqref{eq:endpoint}.

\begin{figure}[t!]
\begin{algorithm}[H]
  \small
  \Fn{extraction}{
  \In{continuity constraint matrices $\mbf{K}^{(i)}$ (size: $\nsum(\nelms) \times (\r_i+1)$) for $i=1,\ldots,\nelms-1$}
  \Out{extraction matrix $\mbf{H}$}
  \bigskip
  $\mbf{H} \gets \mbox{identity matrix}$ (size: $\nsum(\nelms) \times \nsum(\nelms)$)\;
  \For{$i \gets 1$ \KwTo $\nelms-1$}{
    $\mbf{L} \gets \mbf{H}\, \mbf{K}^{(i)}$\;
    \For{$j \gets 0$ \KwTo $\r_i$}{
      $\bar{\mbf{H}} \gets \mbox{sparse null-space of $(j+1)$-th column of } \mbf{L}$\;
      $\mbf{H} \gets \bar{\mbf{H}}\, \mbf{H}$\;
      $\mbf{L} \gets \bar{\mbf{H}}\, \mbf{L}$\;
    }
  }
  }
\end{algorithm}
\caption{Computation of extraction matrix $\mbf{H}$ in pseudo-code.}\label{alg:extraction}
\end{figure}

Using these matrices, we can build the matrix $\mbf{K}^{(i)}$ of size $\nsum(\nelms)\times(\r_i+1)$ which contains all constraints required to enforce $C^{\r_i}$ at $\x_i$. This matrix is defined row-wise in the following manner:
\begin{itemize}
  \item the $(\nsum(i-1)+j)$-th row of $\mbf{K}^{(i)}$ is equal to the $j$-th row of $\mbf{K}^{(i,1)}$;
  \item the $(\nsum(i)+j)$-th row of $\mbf{K}^{(i)}$ is equal to the $j$-th row of $\mbf{K}^{(i,2)}$;
  \item all other rows of $\mbf{K}^{(i)}$ are identically zero.
\end{itemize}
It can be easily verified that for a row vector of coefficients $\mbf{f}$ such that $\mbf{f}\,\mbf{K}^{(i)} = \mbf{0}$, the spline defined by $\mbf{f}\,\mbf{\bs}$ is going to be $C^{\r_i}$ across $\x_i$. Therefore, once all the matrices $\mbf{K}^{(i)}$ have been assembled, the only remaining step is the construction of $\mbf{H}$ such that it spans their left null-spaces. The matrix $\mbf{H}$ is called \emph{multi-degree spline extraction operator}, and we employ the algorithm in Figure~\ref{alg:extraction} for its construction.
%
The algorithm addresses a single continuity constraint at a time, so this increases the smoothness of the basis functions (obtained by \eqref{eq:multiDegBasisExt}) gradually.
The matrix $\mbf{L}$ keeps track of the remaining continuity constraints for the basis functions built so far.

We now focus on the left null-space computation of a column of the continuity constraint matrix $\mbf{L}$.
Any basis of the null-space would lead to a valid basis of the space $\splSpacerp$ using the previously described procedure. However, we are not just interested in any basis, but are looking for the MDTB-spline basis. We employ the algorithm in Figure~\ref{alg:nullspace} for its construction.
It strives to build the sparsest possible left null-space of column vector $\mbf{l}$, containing the next continuity constraint in the matrix $\mbf{L}$. This is equivalent to building the conversion matrix between two MDTB-spline bases of different smoothness as in \eqref{eq:knotremoval}. Hence, the algorithm can follow the pattern described in \eqref{eq:knotremoval-pattern}.

\begin{figure}[t!]
\begin{algorithm}[H]
  \small
  \Fn{null-space}{
  \In{vector $\mbf{l}$ (length $q$)}
  \Out{left null-space matrix $\bar{\mbf{H}}$}
  \bigskip
  $\bar{\mbf{H}} \gets \mbox{zero matrix}$ (size: $(q-1) \times q$)\;
  $i_1 \gets$ index of first non-zero element of $\mbf{l}$\;
  $i_2 \gets$ index of last non-zero element of $\mbf{l}$\;
  \For{$j \gets 1$ \KwTo $i_1$}{
    $\bar{\mbf{H}}(j,j) \gets 1$\;
  }
  \For{$j \gets i_1+1$ \KwTo $i_2-1$}{
    $\bar{\mbf{H}}(j-1,j) \gets -\dfrac{\mbf{l}(j-1)}{\mbf{l}(j)}\,\bar{\mbf{H}}(j-1,j-1)$\;
    $\bar{\mbf{H}}(j,j) \gets 1 - \bar{\mbf{H}}(j-1,j)$\;
  }
  \For{$j \gets i_2$ \KwTo $q$}{
    $\bar{\mbf{H}}(j-1,j) \gets 1$\;
  }
  }
\end{algorithm}
\caption{Computation of sparse null-space $\bar{\mbf{H}}$ of vector $\mbf{l}$ (column of $\mbf{L}$) in pseudo-code.}\label{alg:nullspace}
\end{figure}

\begin{remark}\label{rmk:extraction-matlab}
  The described extraction mechanism can be efficiently encoded by exploiting the sparsity of the involved matrices, similar to the polynomial multi-degree spline case \cite[Remark~7]{Speleers:2019}. Furthermore, imposing periodicity can be easily built into the procedure by circularly shifting the rows of the extraction matrix such that the periodic continuity constraints behave like continuity constraints at an interior segment join \cite[Remark~9]{Speleers:2019}. 
\end{remark}

\subsection{Computation of Bernstein Functions} \label{sec:computation-bernstein}

The only missing aspect in the computation of the extraction matrix $\mbf{H}$ described in Section~\ref{sec:computation-spline} is the explicit construction of the matrices $\mbf{K}^{(i,1)}$ and $\mbf{K}^{(i,2)}$ in \eqref{eq:K_L}--\eqref{eq:K_R}.
This requires the computation of derivatives of \Chef{} Bernstein functions at the end points of their basic interval. Furthermore, once $\mbf{H}$ has been computed, evaluation and manipulation of MDTB-splines boils down to the equivalent operations on \Chef{} Bernstein functions. Hence, we now focus on the practical computation of Bernstein functions. 

Consider the ECT-space $\ECT{\p}{}(\interval)$ on the interval $\interval:=[\x_0,\x_1]$. Suppose we know a set of functions $\{\ts_0,\ldots,\ts_\p\}$ that forms a basis of $\ECT{\p}{}(\interval)$. We arrange these basis functions in a column vector $\mbf{\ts}$ of length $\p+1$. We are looking for the set of Bernstein functions $\{\bs_{0,\p},\ldots,\bs_{\p,\p}\}$ that span the same space, and we arrange these basis functions in a column vector $\mbf{\bs}$ of length $\p+1$. Then, we aim to construct a conversion matrix $\mbf{C}$ of size $(\p+1) \times (\p+1)$ such that 
\begin{equation}\label{eq:conversion}
  \mbf{\bs} = \mbf{C}\,\mbf{\ts}.  
\end{equation}
For the computation of $\mbf{C}$, we can rely on the end-point interpolation properties of Bernstein functions in \eqref{eq:endpoint} and solve a specific Hermite interpolation problem in $\ECT{\p}{}(\interval)$ for each Bernstein function; see Remark~\ref{rmk:Hermite-Bernstein}. Let $\mbf{M}_0$ and $\mbf{M}_1$ be the transposed Wronskian matrix of the basis $\mbf{\ts}$ evaluated at the end points $\x_0$ and $\x_1$, respectively, so
\begin{equation*}
  \mbf{M}_0 := (\mbf{W}(\x_0))^T, 
  \quad
  \mbf{M}_1 := (\mbf{W}(\x_1))^T,
\end{equation*}
where 
\begin{equation*}
  \mbf{W}(\x) := \begin{bmatrix}
 \ts_0(x) & \ts_1(x) & \cdots & \ts_\p(x) \\
 D\ts_0(x) & D\ts_1(x) & \cdots & D\ts_\p(x) \\
 \vdots & \vdots & & \vdots \\
 D^\p\ts_0(x) & D^\p\ts_1(x) & \cdots & D^\p\ts_\p(x) \\ 
\end{bmatrix}.
\end{equation*}
Then, we can employ the algorithm in Figure~\ref{alg:conversion} for the construction of $\mbf{C}$. Note that the use of matrix inverses is just for notational convenience, and they should not be computed explicitly! Instead of inverting a matrix, a linear system should be solved by means of standard numerical linear algebra routines. These linear systems are unisolvent by the properties of the Bernstein functions and ECT-spaces.

\begin{figure}[t!]
\begin{algorithm}[H]
  \small
  \Fn{conversion}{
  \In{end-point derivative matrices $\mbf{M}_0$ and $\mbf{M}_1$ (size: $(\p+1) \times (\p+1)$)}
  \Out{conversion matrix $\mbf{C}$}
  \bigskip
  $\mbf{C} \gets \mbox{zero matrix}$ (size: $(\p+1) \times (\p+1)$)\;
  $\mbf{s} \gets \mbox{zero row vector}$ (size: $1 \times (\p+1)$)\;
  $\mbf{e} \gets \mbox{first canonical row vector}$ (size: $1 \times (\p+1)$)\;
  $\mbf{C}(\p+1,:) \gets \mbf{e}\, [\mbf{M}_1(:, 1), \mbf{M}_0(:, 1\,\mbox{:}\,\p)]^{-1}$\;
  \For{$i \gets 2$ \KwTo $\p$}{
    $\mbf{s} \gets \mbf{s} + \mbf{C}(\p-i+3, :)$\;
    $\mbf{c} \gets \mbox{zero row vector}$ (size: $1 \times (\p+1)$)\;
    $\mbf{c}(i) \gets - \mbf{s}\, \mbf{M}_1(:, i)$\;
    $\mbf{C}(\p-i+2, :) \gets \mbf{c}\,[\mbf{M}_1(:, 1\,\mbox{:}\,i), \mbf{M}_0(:, 1\,\mbox{:}\,(\p-i+1))]^{-1}$\;
  }
  $\mbf{C}(1,:) \gets \mbf{e}\, [\mbf{M}_0(:, 1), \mbf{M}_1(:, 1\,\mbox{:}\,\p)]^{-1}$\;
  }
\end{algorithm}
\caption{Computation of conversion matrix $\mbf{C}$ in pseudo-code.}\label{alg:conversion}
\end{figure}

\begin{remark}
  The choice of the basis $\mbf{\ts}$ in \eqref{eq:conversion} has a major influence on the computation of $\mbf{\bs}$, and hence is of utmost importance. Specific knowledge of the ECT-space, with detection of possible instabilities, and computational efficiency need to be taken into consideration for this choice. Examples are given in Section~\ref{sec:implementation-ECT}.
\end{remark}

\section{Implementational Aspects} \label{sec:implementation}

The construction and use of MDTB-splines through the previously described extraction procedure has been implemented in a small object-oriented \Matlab{} toolbox that is available through the CALGO library. In this section, we give an overview of its class structure.
Full details of the facilities available from the \Matlab{} toolbox may be found in the user manual that accompanies the software. 

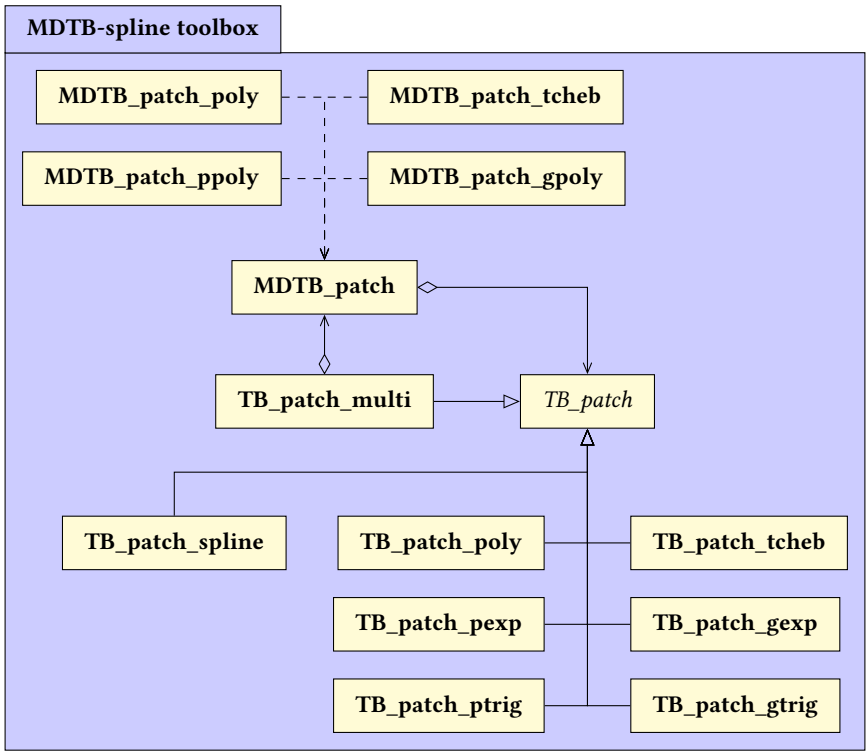
\begin{figure}
  \begin{tikzpicture} 
  \begin{umlpackage}[x=0, y=0, alias=mdtb-spline]{MDTB-spline toolbox}
    \def\s{0.07}
    \umlsimpleclass[x=2-\s, y=8, alias=mdtb-poly, anchor=east]{MDTB\_patch\_poly} 
    \umlsimpleclass[x=2-\s, y=7-\s, alias=mdtb-ppoly, anchor=east]{MDTB\_patch\_ppoly} 
    \umlsimpleclass[x=3+\s, y=8, alias=mdtb-tcheb, anchor=west]{MDTB\_patch\_tcheb}
    \umlsimpleclass[x=3+\s, y=7-\s, alias=mdtb-gpoly, anchor=west]{MDTB\_patch\_gpoly} 
    \umlsimpleclass[x=2.5, y=5.5, alias=mdtb-patch]{MDTB\_patch} 
    \umlsimpleclass[x=6, y=4, alias=tb-patch, type=abstract]{TB\_patch}
    \umlsimpleclass[x=2.5, y=4, alias=tb-multi]{TB\_patch\_multi} 
    \umlsimpleclass[x=0.5, y=2+2*\s, alias=tb-spline]{TB\_patch\_spline} 
    \umlsimpleclass[x=5.5-\s, y=2+2*\s, alias=tb-poly, anchor=east]{TB\_patch\_poly} 
    \umlsimpleclass[x=5.5-\s, y=1+\s, alias=tb-pexp, anchor=east]{TB\_patch\_pexp} 
    \umlsimpleclass[x=5.5-\s, y=0, alias=tb-ptrig, anchor=east]{TB\_patch\_ptrig}
    \umlsimpleclass[x=6.5+\s, y=2+2*\s, alias=tb-tcheb, anchor=west]{TB\_patch\_tcheb} 
    \umlsimpleclass[x=6.5+\s, y=1+\s, alias=tb-gexp, anchor=west]{TB\_patch\_gexp} 
    \umlsimpleclass[x=6.5+\s, y=0, alias=tb-gtrig, anchor=west]{TB\_patch\_gtrig}
    \umlHVdep{mdtb-tcheb}{mdtb-patch}
    \umlHVdep{mdtb-gpoly}{mdtb-patch}
    \umlHVdep{mdtb-poly}{mdtb-patch}
    \umlHVdep{mdtb-ppoly}{mdtb-patch}
    \umlHVuniaggreg{mdtb-patch}{tb-patch}
    \umluniaggreg{tb-multi}{mdtb-patch}
    \umlinherit{tb-multi}{tb-patch}
    \umlVHVinherit{tb-spline}{tb-patch}
    \umlHVinherit{tb-poly}{tb-patch}
    \umlHVinherit{tb-pexp}{tb-patch}
    \umlHVinherit{tb-ptrig}{tb-patch}
    \umlHVinherit{tb-tcheb}{tb-patch}
    \umlHVinherit{tb-gexp}{tb-patch}
    \umlHVinherit{tb-gtrig}{tb-patch}
  \end{umlpackage} 
  \end{tikzpicture}
  \caption{Class diagram of the MDTB-spline toolbox in \Matlab{}.}\label{fig:class-diagram}
\end{figure}

\subsection{Object-Oriented Implementation} \label{sec:implementation-oo}

The class diagram of the MDTB-spline toolbox is shown in Figure~\ref{fig:class-diagram}. The central class in the toolbox is the class \texttt{MDTB\_patch}, which provides functionality for computing the MDTB-spline extraction matrix, both in the periodic and non-periodic spline setting. Furthermore, it allows for evaluating, differentiating and visualizing the obtained MDTB-spline basis functions and any MDT-spline function represented in such basis.

The class \texttt{MDTB\_patch} is built upon the class \texttt{TB\_patch}, which mainly identifies a local ECT-space. According to the definition of MDT-spline spaces (see Definition~\ref{def:spline-space}), each object of type \texttt{MDTB\_patch} contains a (heterogeneous) array of objects of type \texttt{TB\_patch}. In order to reflect the heterogeneous nature of ECT-spaces, the class \texttt{TB\_patch} is \emph{abstract}. The functionality of evaluation and differentiation of the Bernstein basis functions is delegated to (specialized) child classes, as well as the option to provide a separate implementation for end-point derivatives because of their importance in the MDTB-spline framework.
Manipulation of functions represented in the Bernstein basis and their visualization is handled in the class \texttt{TB\_patch} itself.

As indicated in the class diagram in Figure~\ref{fig:class-diagram}, there are several child classes of the class \texttt{TB\_patch} available in the toolbox. They provide functionality to work with general ECT-spaces based on constant-coefficient linear differential operators (class \texttt{TB\_patch\_tcheb}), algebraic polynomial spaces (classes \texttt{TB\_patch\_poly} and \texttt{TB\_patch\_spline}), other polynomial-type spaces (classes \texttt{TB\_patch\_pexp} and \texttt{TB\_patch\_ptrig}), and generalized polynomial spaces (classes \texttt{TB\_patch\_gexp} and \texttt{TB\_patch\_gtrig}). In Section~\ref{sec:implementation-ECT} we discuss these example spaces and their implementations in more detail. 
Thanks to the object-oriented structure of the toolbox, other ECT-spaces and/or specialized implementations can be easily incorporated by adding new child classes of \texttt{TB\_patch} --- they need to implement the three methods \texttt{TB\_evaluation\_all}, \texttt{TB\_differentiation\_all}, and \texttt{TB\_diffend\_all}).

The purpose of the class \texttt{TB\_patch\_multi} is to encapsulate an object of the class \texttt{MDTB\_patch} and a non-periodic extraction matrix such that the corresponding multi-degree spline space can be treated as if it is an instance of type \texttt{TB\_patch}. In this way, already constructed multi-degree spline spaces can be embedded into larger multi-degree spline spaces without the need for recomputing them.

The classes \texttt{MDTB\_patch\_tcheb}, \texttt{MDTB\_patch\_poly}, \texttt{MDTB\_patch\_ppoly}, and \texttt{MDTB\_patch\_gpoly} are factory classes for \texttt{MDTB\_patch}. They provide simplified functionality to initialize objects of type \texttt{MDTB\_patch} consisting of local ECT-spaces based on constant-coefficient linear differential operators, algebraic polynomial spaces, other polynomial-type spaces, and generalized polynomial spaces, respectively.

\begin{remark}
  This MDTB-spline toolbox is a redesigned, object-oriented and extended version of the MDB-spline toolbox provided by \citet{Speleers:2019}, also implemented in \Matlab{}. The latter toolbox only deals with polynomial MDB-splines, and its functionality is essentially covered by the classes \texttt{TB\_patch}, \texttt{TB\_patch\_spline}, \texttt{MDTB\_patch} and \texttt{MDTB\_patch\_poly}. Even though the MDB-spline toolbox is not object-oriented, the syntax of the new MDTB-spline toolbox is (almost fully) compatible with it and all its function calls are still available (under the minor restriction that the names start now with ``\texttt{(MD)TB\_*}'' instead of ``\texttt{(MD)B\_*}'', so as to emphasize the \Chef{} nature of the extension). 
\end{remark}

\subsection{Implementation of ECT-Spaces} \label{sec:implementation-ECT}

Here we focus on the large class of ECT-spaces described in Section~\ref{sec:ECT-Np}.
We refer the reader to \citet{Roth:2019} for its practical relevance.
Consider the null-space $\NN_\p$ of a constant-coefficient linear differential operator $\cL_\p$ as in \eqref{eq:Lp}  with $\p\geq1$. We assume that $\NN_\p$ is an ECT-space (so the interval $\interval:=[\x_0,\x_1]$ is chosen sufficiently small) and that $\croot=0$ is at least a first-order root of the characteristic polynomial \eqref{eq:Lp-pol}. This ECT-space can be uniquely identified by means of the triples
\begin{equation}\label{eq:roots}
  (0,0,\mu_0), \quad (\alpha_1,\beta_1,\mu_1), \quad \ldots, \quad (\alpha_m,\beta_m,\mu_m),
\end{equation}
representing the different roots $\croot_k=\alpha_k+\ii \beta_k$ of order $\mu_k$ ($\mu_k\geq1$) of the polynomial \eqref{eq:Lp-pol}. 
Such space is translation-invariant. When we want to emphasize the specific parameters in \eqref{eq:roots}, the space will be denoted with $\NN_\p^{(\alpha_1,\beta_1,\mu_1), \ldots, (\alpha_m,\beta_m,\mu_m)}$. For simplicity of notation, we assume that the complex conjugate roots are excluded from \eqref{eq:roots}.

The basis $\mbf{\ts}$ in \eqref{eq:conversion} depends on the type of roots in \eqref{eq:roots} and is chosen as follows. We differentiate between four types of roots for implementational efficiency. A given root $\croot=\alpha+\ii \beta$ of order $\mu$ gives rise to the following basis functions for $i=0,\ldots,\mu-1$:
\begin{itemize}
  \item if $\alpha=0$ and $\beta=0$, then
  \begin{equation*}
   \ts_i(\x) = \frac{(\x-\x_0)^i}{i!}, \quad 
   D\ts_i(\x) = \ts_{i-1}(\x);
  \end{equation*}
  \item if $\alpha\neq0$ and $\beta=0$, then
  \begin{equation*}
   \ts_i(\x) = \frac{(\x-\x_0)^i}{i!}\ee^{\alpha(\x-\x_0)}, \quad 
   D\ts_i(\x) = \ts_{i-1}(\x)+\alpha\ts_{i}(\x);
  \end{equation*}
  \item if $\alpha=0$ and $\beta\neq0$, then
  \begin{equation*}
   \ts_{2i}(\x) = \frac{(\x-\x_0)^i}{i!}\cos(\beta(\x-\x_0)), \quad 
   \ts_{2i+1}(\x) = \frac{(\x-\x_0)^i}{i!}\sin(\beta(\x-\x_0)),
  \end{equation*}
  and their derivatives
  \begin{equation*}
    D\ts_{2i}(\x) = \ts_{2i-2}(\x)-\beta\ts_{2i+1}(\x), \quad
    D\ts_{2i+1}(\x) = \ts_{2i-1}(\x)+\beta\ts_{2i}(\x);
  \end{equation*}
  \item if $\alpha\neq0$ and $\beta\neq0$, then
  \begin{equation*}
   \ \ \ \ts_{2i}(\x) = \frac{(\x-\x_0)^i}{i!}\ee^{\alpha(\x-\x_0)}\cos(\beta(\x-\x_0)), 
   \quad \ts_{2i+1}(\x) = \frac{(\x-\x_0)^i}{i!}\ee^{\alpha(\x-\x_0)}\sin(\beta(\x-\x_0)),
  \end{equation*}
  and their derivatives
  \begin{equation*}
  \ \ \ D\ts_{2i}(\x) = \ts_{2i-2}(\x)+\alpha\ts_{2i}(\x)-\beta\ts_{2i+1}(\x), \quad
    D\ts_{2i+1}(\x) = \ts_{2i-1}(\x)+\alpha\ts_{2i+1}(\x)+\beta\ts_{2i}(\x).
  \end{equation*}
\end{itemize}
In the expressions of the derivatives it is assumed that $\ts_i(\x)=0$ for $i<0$.
These basis functions are chosen so that the conversion matrix $\mbf{C}$ will be invariant under translations of the interval~$\interval$. Moreover, when evaluating all basis functions (and their derivatives) simultaneously, a lot of arithmetic calculations can be shared due to their particular structure.

This class of ECT-spaces is extremely flexible, as there are plenty of (shape) parameters to be chosen in \eqref{eq:roots}. However, this flexibility comes at a risk because it is very easy to select a combination of shape parameters that leads to numerical instabilities, e.g., when choosing two roots that are different but very close to each other. Such instabilities will lead to linear systems that are ill-conditioned, and hence their numerical solution might not be accurate anymore; see Section~\ref{sec:instability} for a discussion about instabilities.
In the \Matlab{} toolbox, the class \texttt{TB\_patch\_tcheb} deals with such general ECT-spaces.

In the following three subsections, we discuss subclasses of ECT-spaces of particular interest, for which an improved implementation has been provided that is more efficient and/or more robust. In the \Matlab{} toolbox, they have been addressed in separate \Matlab{} classes (see the class diagram in Figure~\ref{fig:class-diagram}) and should be selected whenever possible.

\begin{remark} \label{rmk:Roth}
  The same Bernstein functionality as in the \Matlab{} class \texttt{TB\_patch\_tcheb} has also been implemented in the C++ library by \citet{Roth:2019}. However, the implementation of the class \texttt{TB\_patch\_tcheb} differs at the following points:
  \begin{itemize}
    \item The Bernstein basis is directly computed from the basis $\mbf{\ts}$ in \eqref{eq:conversion}, instead of first building an intermediate bicanonical basis and then obtaining the Bernstein functions via a particular LU-decomposition of the corresponding Wronskian matrix. 
    \item The (high-order) derivatives of the basis functions in $\mbf{\ts}$ are computed by means of the above recurrence relations instead of application of the general Leibniz rule.
    \item Since all Bernstein basis functions are treated simultaneously (for their construction, evaluation, and differentiation), a lot of arithmetic calculations can be shared. Moreover, all parts of the \Matlab{} code are highly vectorized. Finally, because of their importance in the MDTB-spline framework, end-point derivatives are implemented separately, resulting in additional computational speed and stability.
  \end{itemize}
\end{remark}

\subsubsection{Algebraic Polynomial Spaces} \label{sec:implementation-ECT-poly}
Algebraic polynomial spaces are the most established 
ECT-spaces (see Examples~\ref{ex:poly} and \ref{ex:nullspace-poly}). In this case, the \Chef{} Bernstein functions are the classical Bernstein polynomials (see Example~\ref{ex:poly-bernstein}). Since we know a simple explicit expression for them, it is not necessary to use the conversion procedure described in Section~\ref{sec:computation-bernstein}. Alternatively, they can be evaluated at $x\in[\x_0,\x_1]$ through the following stable recurrence relation:
\begin{equation*}
  \bs_{j,\p}(\x)=\left(\frac{\x-\x_0}{\x_1-\x_0}\right)\bs_{j-1,\p-1}(\x)+\left(\frac{\x_1-\x}{\x_1-\x_0}\right)\bs_{j,\p-1}(\x),
\quad \p\geq1,
\end{equation*}
with $\bs_{j,\p}(\x)=0$ for $j<0$, $j>\p$ and $\bs_{0,0}(\x)=1$.
Their derivatives can be easily computed as 
\begin{equation*}
  D\bs_{j,\p}(\x)=\frac{p}{\x_1-\x_0}\bigl(\bs_{j-1,\p-1}(\x)-\bs_{j,\p-1}(\x)\bigr).
\end{equation*}
This has been implemented in the \Matlab{} class \texttt{TB\_patch\_poly}.

When dealing with spline spaces consisting of local algebraic polynomial spaces of the same degree $p$, then the MDTB-splines are nothing but the classical polynomial B-splines. In this case, instead of applying the extraction process described in Section~\ref{sec:computation-spline}, one can also use the following stable recurrence relation for their computation. Given an \emph{open knot vector} of the form
\begin{equation} \label{eq:knots}
\mbf{\xi}{} := (\xi_k )_{k=1}^{\N+\p+1}:= (\;
    \underbrace{\x_0,\; \ldots,\;  \x_0}_{\p - \r_0 \text{ times}},\;  \ldots,\;
    \underbrace{\x_{i},\;  \ldots,\;  \x_{i}}_{\p - \r_i \text{ times}},\;   \ldots,\;
    \underbrace{\x_{\nelms},\;  \ldots,\;  \x_{\nelms}}_{\p - \r_{\nelms} \text{ times}}),
\end{equation}
the polynomial B-splines $\BS_{k,\p}$, $k=1,\ldots,\N$ can be evaluated at $\x \in [\x_0,\x_{\nelms})$ as 
\begin{equation*}
  \BS_{k,\p}(\x) = \frac{\x-\xi_k}{\xi_{k+\p}-\xi_k}\,\BS_{k,\p-1}(\x) + \frac{\xi_{k+\p+1}-\x}{\xi_{k+\p+1}-\xi_{k+1}}\,\BS_{k+1,\p-1}(\x),
\end{equation*}
starting from
\begin{equation*}
  \BS_{k,0}(\x) = 
  \begin{cases}
    1, & \text{if }\xi_k \leq \x < \xi_{k+1},\\
    0, & \text{otherwise},
  \end{cases}
\end{equation*}
and under the convention that fractions with zero denominator have value zero. 
Their derivatives can be computed as
\begin{equation*}
  D_+\BS_{k,\p}(\x) = p\left(\frac{\BS_{k,\p-1}(\x)}{\xi_{k+\p}-\xi_k}-\frac{\BS_{k+1,\p-1}(\x)}{\xi_{k+\p+1}-\xi_{k+1}}\right).
\end{equation*}
Such spline space can be used as local space in the MDTB-spline setting, and will be more efficient than working with separate polynomial spaces of the same degree \cite[Remark~3]{Speleers:2019}. This has been implemented in the \Matlab{} class \texttt{TB\_patch\_spline}; its implementation has been borrowed from the MDB-spline toolbox developed by \citet{Speleers:2019}.

\begin{remark}
The MDTB-spline knot vectors in \eqref{eq:knots-u}--\eqref{eq:knots-v} relate to the classical B-spline knot vector in \eqref{eq:knots} as $\lknot_k = \knot_k$ and $\rknot_k = \knot_{k+\p+1}$, $k = 1, \ldots, \N$.
\end{remark}

\subsubsection{Other ECT-spaces with stable recurrence relations} \label{sec:implementation-ECT-ppoly}

Besides algebraic polynomial spaces, there are few other classes of ECT-spaces for which simple and stable recurrence relations are known for the evaluation of the \Chef{} Bernstein functions and their derivatives.
In particular, in the \Matlab{} toolbox we consider two classes of polynomial-type spaces of the form
\begin{equation*}
\PP_\p^{\mathfrak{u},\mathfrak{v}}
:= \bigl\langle 1,\mathfrak{u}(\x),\mathfrak{v}(\x),\mathfrak{u}(2\x),\mathfrak{v}(2\x),\ldots,\mathfrak{u}(q\x),\mathfrak{v}(q\x)\bigr\rangle, \quad p=2q\geq2,
\end{equation*}
where $\mathfrak{u},\mathfrak{v}$ are chosen to be exponential or trigonometric functions.
The first type of space corresponds to the null-space in Example~\ref{ex:nullspace-exp} with $\p=m=2q$, $\mu_k=1$, $k=0,\ldots,2q$, and $\alpha_{2k-1}=-\alpha_{2k}=k\alpha$, $k=1,\ldots,q$, and we denote it with
\begin{align*}
\PE_\p^\alpha &:= \NN_\p^{(\alpha,0,1), (-\alpha,0,1),\ldots, (q\alpha,0,1), (-q\alpha,0,1)} 
=\bigl\langle 1,\ee^{\alpha\x},\ee^{-\alpha\x},\ee^{2\alpha\x},\ee^{-2\alpha\x},\ldots,\ee^{q\alpha\x},\ee^{-q\alpha\x}\bigr\rangle\\
&\;=\bigl\langle 1,\cosh(\alpha\x),\sinh(\alpha\x),\cosh(2\alpha\x),\sinh(2\alpha\x),\ldots,\cosh(q\alpha\x),\sinh(q\alpha\x)\bigr\rangle.
\end{align*}
The second type of space is the null-space 
\begin{equation*}
\PT_\p^\beta := \NN_\p^{(0,\beta,1),\ldots, (0,q\beta,1)}
=\bigl\langle 1,\cos(\beta\x),\sin(\beta\x),\cos(2\beta\x),\sin(2\beta\x),\ldots,\cos(q\beta\x),\sin(q\beta\x)\bigr\rangle.
\end{equation*}
They are both invariant under translations and reflections. 

To address these two spaces, we set
\begin{equation*}
  \mathfrak{u}(\x)=\cosh(\alpha\x), \quad \mathfrak{v}(\x)=\sinh(\alpha\x)
\end{equation*}
and
\begin{equation*}
  \mathfrak{u}(\x)=\cos(\beta\x), \quad \mathfrak{v}(\x)=\sin(\beta\x),
\end{equation*}
respectively.
Then, we can evaluate the \Chef{} Bernstein functions at 
$x\in[\x_0,\x_1]$ through the following stable recurrence relation:
\begin{equation*}
  \bs_{j,\p}(\x)=\bs_{2,2}(\x)\,\bs_{j-2,\p-2}(\x)+\bs_{1,2}(\x)\,\bs_{j-1,\p-2}(\x)+\bs_{0,2}(\x)\,\bs_{j,\p-2}(\x), \quad \p\geq4,
\end{equation*}
with $\bs_{j,\p}(\x)=0$ for $j<0$, $j>\p$. Moreover, we have
\begin{equation*}
  \bs_{0,2}(\x) = \sigma_{0,2}\, \mathfrak{v}^2\left(\frac{\x_1-\x}{2}\right), \quad
  \bs_{1,2}(\x) = \sigma_{1,2}\, \mathfrak{v}\left(\frac{\x_1-\x}{2}\right) \mathfrak{v}\left(\frac{\x-\x_0}{2}\right), \quad
  \bs_{2,2}(\x) = \sigma_{2,2}\, \mathfrak{v}^2\left(\frac{\x-\x_0}{2}\right),
\end{equation*}
and
\begin{equation*}
  \sigma_{0,2}=\sigma_{2,2}= 1\Big / \mathfrak{v}^2\left(\frac{\x_1-\x_0}{2}\right), \quad
  \sigma_{1,2}= 2\,\mathfrak{u}\left(\frac{\x_1-\x_0}{2}\right)\Big / \mathfrak{v}^2\left(\frac{\x_1-\x_0}{2}\right).
\end{equation*}
Their derivatives can be computed as 
\begin{equation*}
  D\bs_{j,\p}(\x) = \tau\left(j \frac{\sigma_{2,2}\sigma_{j,p}}{\sigma_{j-1,p}}\,\bs_{j-1,\p}(\x)
  -(q-j) \sigma_{1,2}\,\bs_{j,\p}(\x)
  -(p-j) \frac{\sigma_{0,2}\sigma_{j,p}}{\sigma_{j+1,p}}\,\bs_{j+1,\p}(\x)\right),
\end{equation*}
where
\begin{equation*}
  \tau=\frac{D\mathfrak{v}(0)}{2}\,\mathfrak{v}\left(\frac{\x_1-\x_0}{2}\right),
\end{equation*}
and
\begin{equation*}
  \sigma_{j,\p} = \sigma_{2,2}\,\sigma_{j-2,\p-2}+\sigma_{1,2}\,\sigma_{j-1,\p-2}+\sigma_{0,2}\,\sigma_{j,\p-2},\quad \p\geq4,
\end{equation*}
with $\sigma_{j,\p}=0$ for $j<0$, $j>\p$. More details can be found in \citet{Sanchez:1998} and \citet{Shen:2005}.
These recurrence relations have been implemented in the \Matlab{} classes \texttt{TB\_patch\_pexp} and \texttt{TB\_patch\_ptrig} for the exponential and trigonometric cases, respectively.

\subsubsection{Generalized Polynomial Spaces} \label{sec:implementation-ECT-gpoly}
Generalized polynomial spaces are an important class of ECT-spaces. They can be seen as the minimal extension of algebraic polynomial spaces still offering a wide variety of flexibility. They are defined as
\begin{equation}\label{eq:nullspace-gb}
\GG_\p^{\mathfrak{u},\mathfrak{v}}
:= \bigl\langle 1,x,\ldots,x^{\p-2},\mathfrak{u}(\x),\mathfrak{v}(\x)\bigr\rangle, \quad \p\geq2, 
\end{equation}
for given functions $\mathfrak{u},\mathfrak{v}\in C^\p(\interval)$ such that $\langle D^{\p-1}\mathfrak{u},D^{\p-1}\mathfrak{v}\rangle$ is an ECT-space on $\interval:=[\x_0,\x_1]$.
We refer the reader to \citet{Costantini:2005} and \citet{Lyche:2019} for more details. 
Two practically relevant instances of \eqref{eq:nullspace-gb} are obtained by choosing $\mathfrak{u},\mathfrak{v}$ to be exponential or trigonometric functions. They are implemented in the \Matlab{} classes \texttt{TB\_patch\_gexp} and \texttt{TB\_patch\_gtrig}, respectively. 

First, we consider the null-space in Example~\ref{ex:nullspace-exp} 
with $\p\geq2$, $m=2$, $\mu_0=\p-1$, $\mu_1=\mu_2=1$, and $\alpha_1=-\alpha_2=\alpha$, and denote it with 
\begin{equation*}
\EE_\p^\alpha := \NN_\p^{(\alpha,0,1), (-\alpha,0,1)}
=\bigl\langle 1,x,\ldots,x^{\p-2},\ee^{\alpha\x},\ee^{-\alpha\x}\bigr\rangle 
= \bigl\langle 1,x,\ldots,x^{\p-2},\cosh(\alpha\x),\sinh(\alpha\x)\bigr\rangle.  
\end{equation*}
This is the algebraic polynomial space of degree $\p-2$, enriched with two exponential functions. It is invariant under translations and reflections. For values of $\alpha$ away from zero, we take the basis $\mbf{\ts}$ in \eqref{eq:conversion} as
\begin{equation*}
  \ts_i(\x) = \alpha^i\frac{(\x-\x_0)^i}{i!}, \quad i=0,\ldots,\p-2, 
  \quad
  \ts_{\p-1}(\x) = \ee^{\alpha(\x-\x_0)},
  \quad
  \ts_{\p}(\x) = \ee^{-\alpha(\x-\x_0)}.
\end{equation*}
Note that $D\ts_0=0$, $D\ts_i=\alpha\ts_{i-1}$, $i=1,\ldots,\p-2$, $D\ts_{\p-1}=\alpha\ts_{\p-1}$ and $D\ts_{\p}=-\alpha\ts_{\p}$. This allows us to scale the $i$-th column of both matrices $\mbf{M}_0$ and $\mbf{M}_1$ in Figure~\ref{alg:conversion} by the factor $\alpha^{-i}$ to balance better the overall magnitude of the entries when solving the involved linear systems.
Unfortunately, the above choice of basis becomes numerically unstable when $\alpha$ is close to zero. Therefore, taking into account that $\lim_{\alpha\rightarrow0}\EE_\p^\alpha=\PP_\p$, we employ the following basis for values of $\alpha$ close to zero:
\begin{equation*}
  \ts_i(\x) = \frac{(\x-\x_0)^i}{i!}, \quad i=0,\ldots,\p-2,
  \quad
  \ts_{i}(\x) = \sum_{j=0}^{+\infty}\alpha^{2j}\frac{(\x-\x_0)^{i+2j}}{({i+2j})!},
  \quad
  i=\p-1,\p.
\end{equation*}
The choice of $\ts_{\p-1}$ and $\ts_{\p}$ is based on the Taylor expansion of $\cosh(\alpha(\x-\x_0))$ and $\sinh(\alpha(\x-\x_0))$ around $x=\x_0$. Of course, in practice, we need to truncate these series; the order of truncation should be based on the required level of accuracy and the value of $\alpha$. This dual choice of basis has been implemented in the \Matlab{} class \texttt{TB\_patch\_gexp}.

Second, we consider the null-space in Example~\ref{ex:nullspace-trig} 
with $\p\geq2$, and denote it with 
\begin{equation*}
\TT_\p^\beta:=\NN_\p^{(0,\beta,1)}
=\bigl\langle 1,x,\ldots,x^{\p-2},\cos(\beta\x),\sin(\beta\x)\bigr\rangle.
\end{equation*}
This is the algebraic polynomial space of degree $\p-2$, enriched with two trigonometric functions. It is invariant under translations and reflections. The implementation of the corresponding Bernstein basis functions can be done in a manner that is very similar to the exponential case. For values of $\beta$ away from zero, we take
\begin{equation*}
  \ts_i(\x) = \beta^i\frac{(\x-\x_0)^i}{i!}, \quad i=0,\ldots,\p-2, 
  \quad
  \ts_{\p-1}(\x) = \cos(\beta(\x-\x_0)),
  \quad
  \ts_{\p}(\x) = \sin(\beta(\x-\x_0)),
\end{equation*}
while for values of $\beta$ close to zero, we take
\begin{equation*}
  \ts_i(\x) = \frac{(\x-\x_0)^i}{i!}, \quad i=0,\ldots,\p-2,
  \quad
  \ts_{i}(\x) = \sum_{j=0}^{+\infty}(-1)^j\alpha^{2j}\frac{(\x-\x_0)^{i+2j}}{({i+2j})!},
  \quad
  i=\p-1,\p.
\end{equation*}
This dual choice of basis has been implemented in the \Matlab{} class \texttt{TB\_patch\_gtrig}.

\section{Numerical Examples} \label{sec:examples}

In this section, we illustrate the \Matlab{} toolbox with several examples. We start with the main purpose of the toolbox: computation and manipulation of MDTB-splines. Then, we discuss some of the numerical pitfalls of working with ECT-spaces. We end with a simple application to support theoretical analysis of ECT-spaces and MDT-spline spaces.

\subsection{Some Sets of Basis Functions}

The main purpose of the \Matlab{} toolbox is to provide a flexible way for computing and manipulating MDTB-splines with pieces drawn from different ECT-spaces. The class structure described in Section~\ref{sec:implementation-oo} enables access to both generic implementations and specialized ones fine-tuned for certain ECT-spaces. Here we restrict ourselves to two small examples illustrating this flexibility. Other examples may be found in \citet[Section~8]{Hiemstra:2020} and in the user manual that accompanies the \Matlab{} toolbox. 
We also refer the reader to \citet[Section~4]{Speleers:2019} for specific examples related to polynomial MDB-splines.

We start by showing how easy it is to smoothly join different types of ECT-spaces 
and to compute the corresponding MDTB-spline basis functions. 

\begin{figure}[t!]
  \subfigure[Local \Chef{} Bernstein bases]
  {\includegraphics[height=2.5cm]{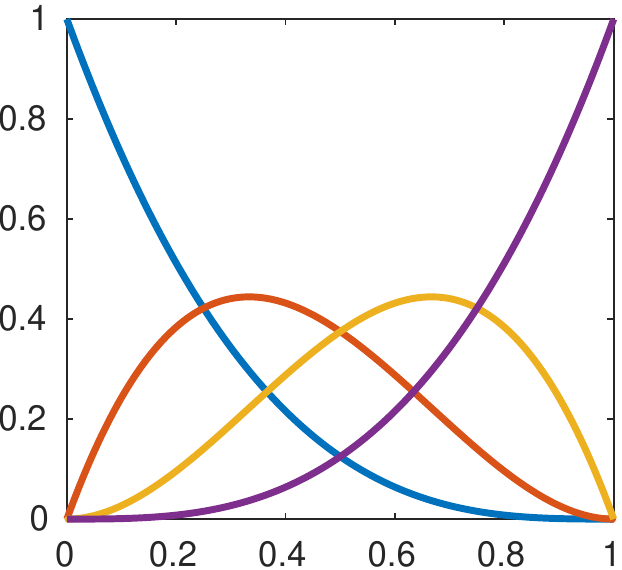} \hspace*{0.1cm}
   \includegraphics[height=2.5cm]{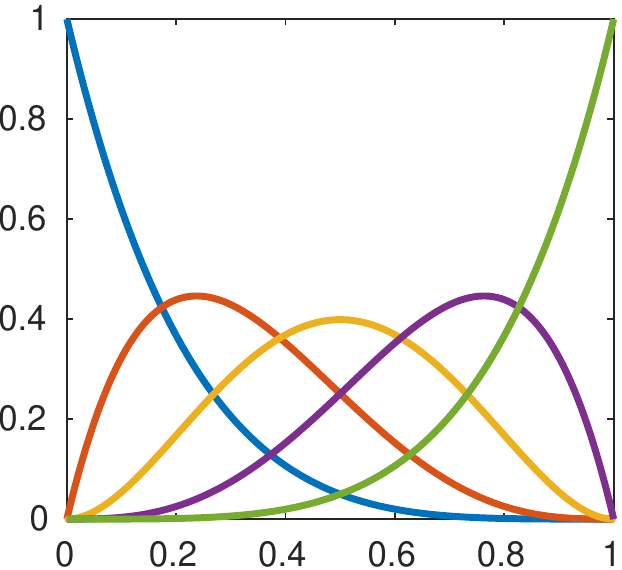} \hspace*{0.1cm}
   \includegraphics[height=2.5cm]{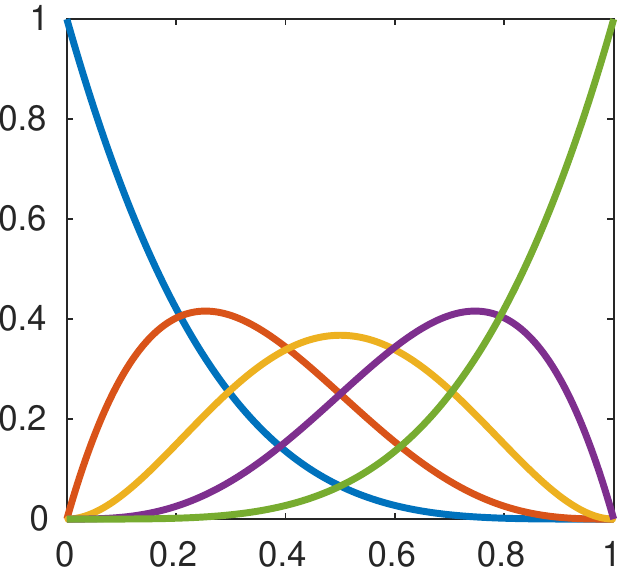} \hspace*{0.1cm}
   \includegraphics[height=2.5cm]{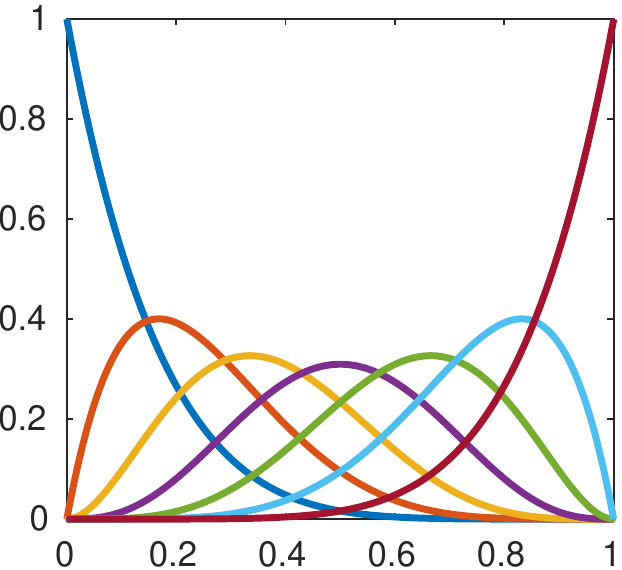}} \\
  \subfigure[Non-periodic MDTB-splines]
  {\includegraphics[height=4.4cm]{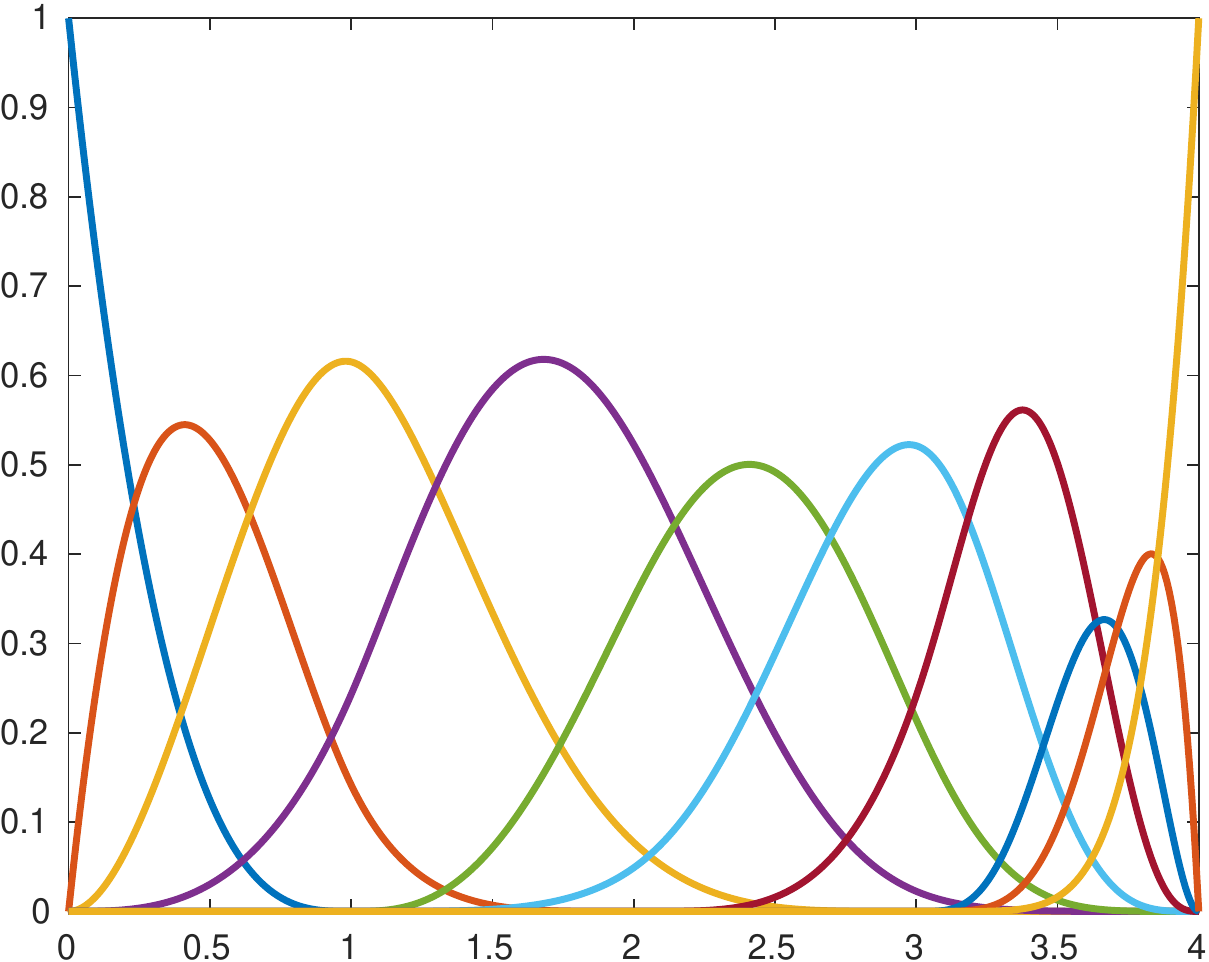}} \hspace*{0.5cm}
  \subfigure[Periodic MDTB-splines]
  {\includegraphics[height=4.4cm]{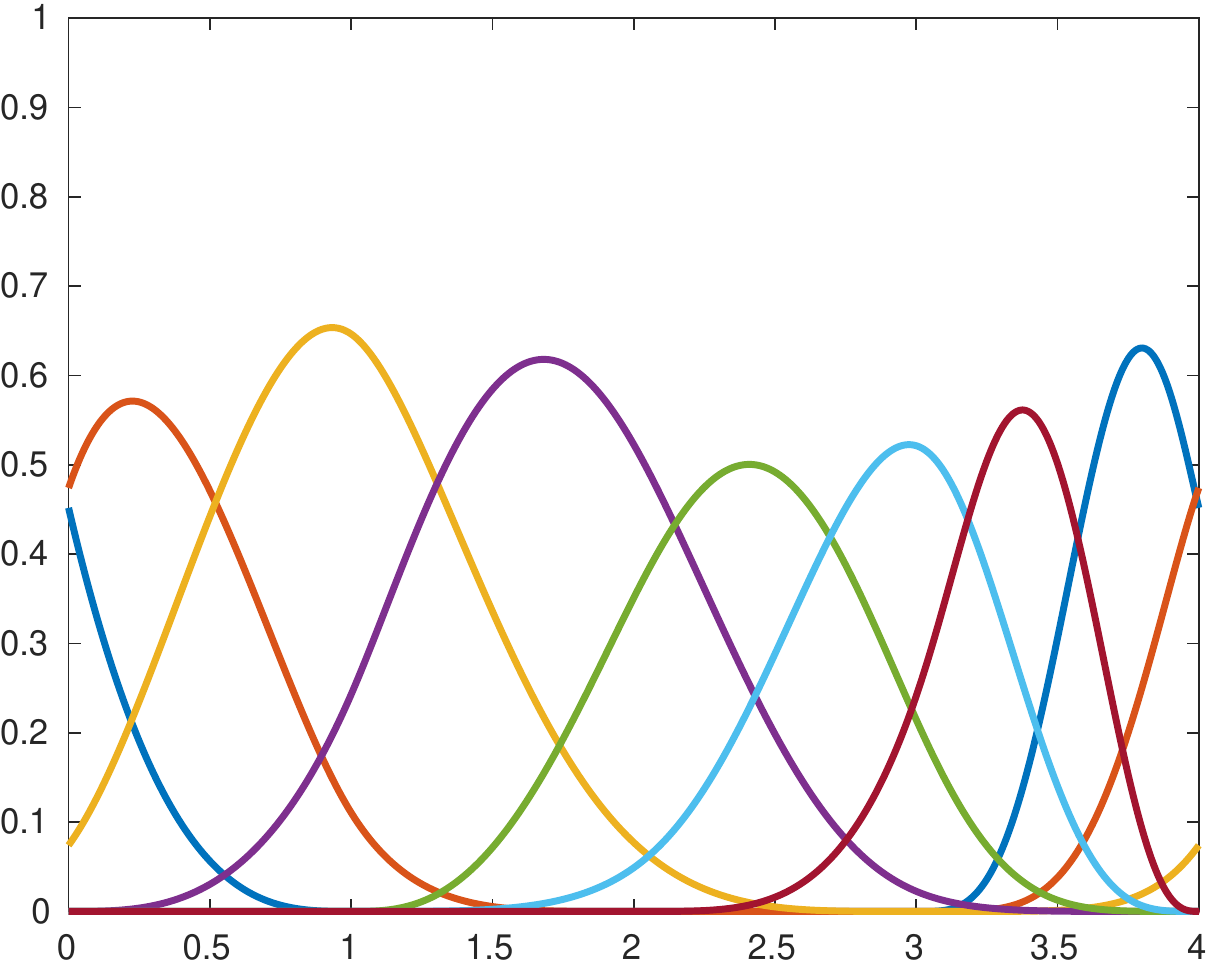}}
  \caption{Non-periodic and periodic MDTB-spline basis functions constructed from local sets of \Chef{} Bernstein basis functions according to Example~\ref{ex:example-basis}.}\label{fig:example-basis}
\end{figure}

\begin{example} \label{ex:example-basis}
  Consider the MDT-spline space $\splSpacerp$ defined by
  \begin{equation*}
    \domain = \{0, 1, 2, 3, 4\}, \quad
    \bp = \{3, 4, 4, 6\}, \quad
    \br = \{-1,2,3,3,-1\},
  \end{equation*}
  and the local ECT-spaces
  \begin{equation*}
    \ECT{3}{(1)} = \PP_3, \quad
    \ECT{4}{(2)} = \EE_4^{3}, \quad
    \ECT{4}{(3)} = \TT_4^{3/2}, \quad
    \ECT{6}{(4)} = \NN_6^{(1,0,1),(-1,0,1),(0,2,1)}.
  \end{equation*}
  These local spaces can be represented as objects of the \Matlab{} classes \texttt{TB\_patch\_poly}, \texttt{TB\_patch\_gexp}, \texttt{TB\_patch\_gtrig} and \texttt{TB\_patch\_tcheb}, respectively. The corresponding \Chef{} Bernstein basis functions are visualized in Figure~\ref{fig:example-basis}(a).
  Imposing the required smoothness at the break points as described by $\br$ gives rise to the MDTB-spline basis functions depicted in Figure~\ref{fig:example-basis}(b).
  The \Matlab{} toolbox also supports periodic MDTB-splines; see Remark~\ref{rmk:extraction-matlab}. Imposing an additional periodic smoothness of $\r_{{\rm per}}=2$ results in the basis functions depicted in Figure~\ref{fig:example-basis}(c). The four central basis functions remain the same as they trivially satisfy the periodic smoothness --- their derivatives up to order two are zero at the end points of the interval.
  This example can be reproduced by executing the \Matlab{} script \texttt{EX\_basis\_A.m}.
\end{example}

MDTB-spline representations are particularly suited for geometric modeling.
The availability of trigonometric and exponential functions allows for the exact description of conic section segments, which can be smoothly blended with polynomial segments. These are very common profiles in industrial design. In the following example, we provide a $C^1$ smooth description of a unit square with adaptable rounded corners (represented as arcs of circles).

\begin{figure}[t!]
  \subfigure[Periodic MDTB-splines for $s=-2,-1,0,1$]
  {\includegraphics[height=2.5cm]{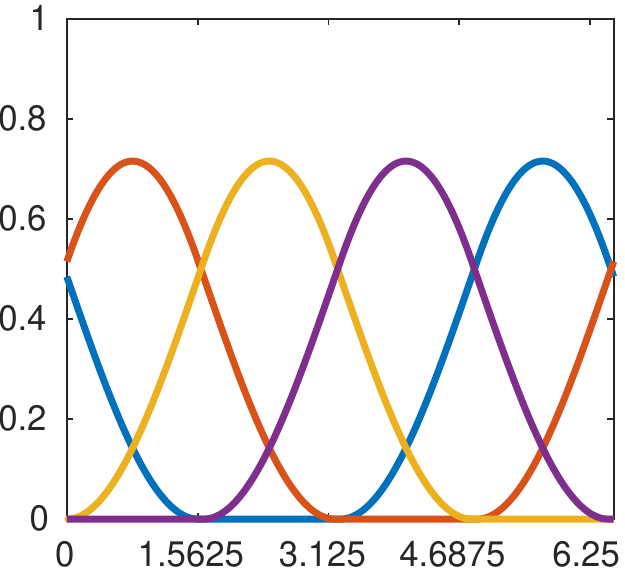} \hspace*{0.1cm}
   \includegraphics[height=2.5cm]{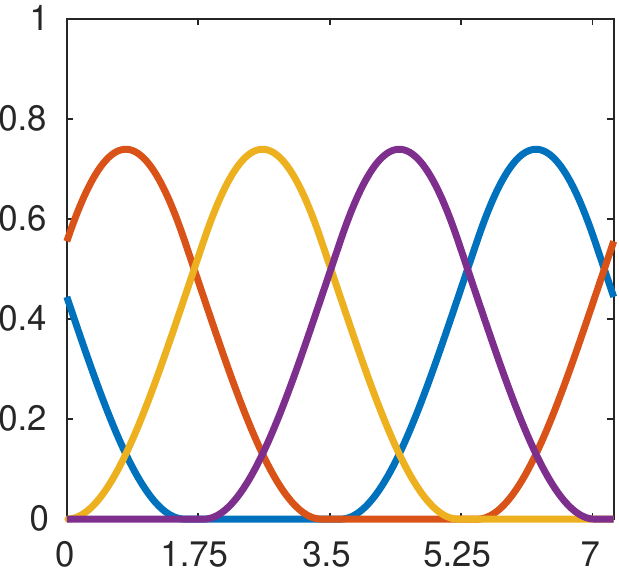} \hspace*{0.1cm}
   \includegraphics[height=2.5cm]{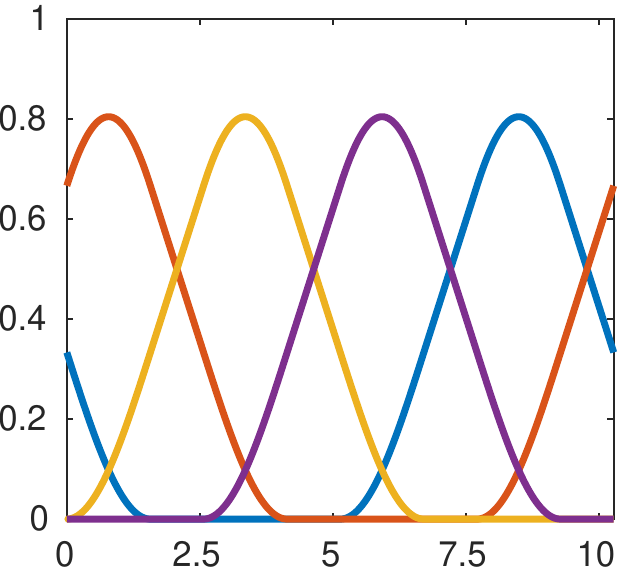} \hspace*{0.1cm}
   \includegraphics[height=2.5cm]{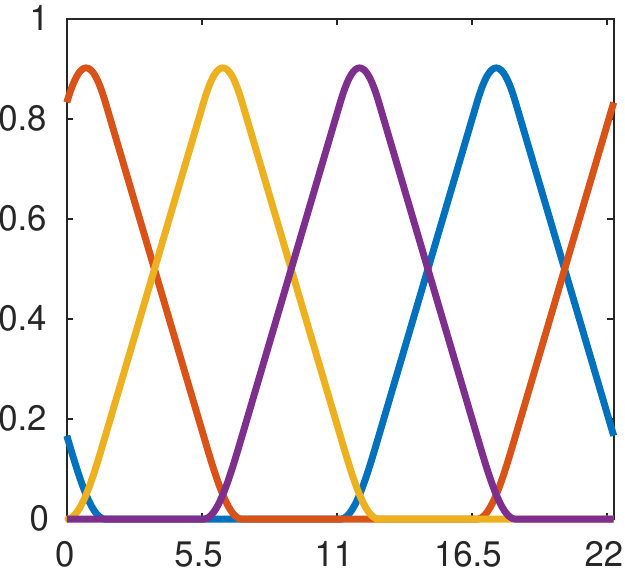}} \\
  \subfigure[Square with rounded corners for $s=-2,-1,0,1$]
  {\includegraphics[height=2.5cm]{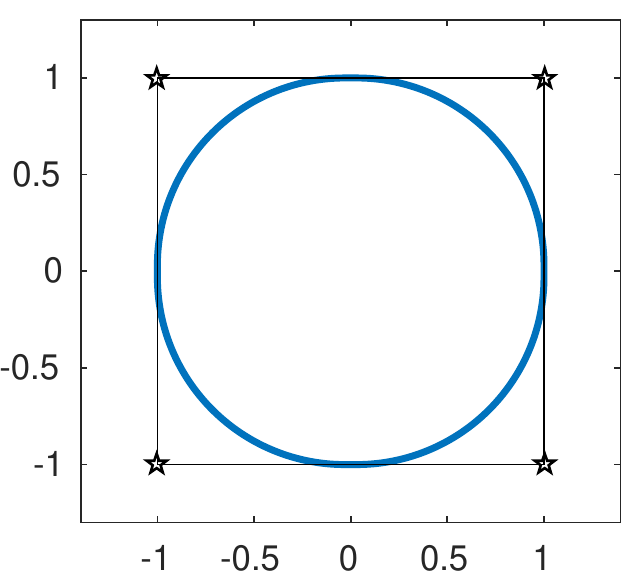} \hspace*{0.1cm}
   \includegraphics[height=2.5cm]{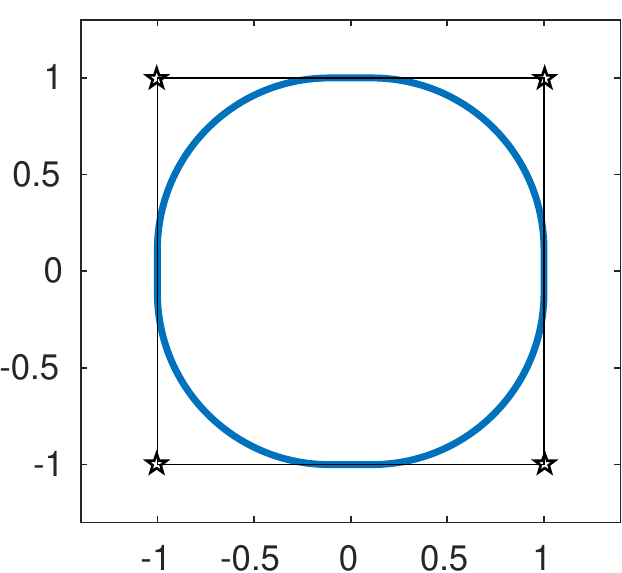} \hspace*{0.1cm}
   \includegraphics[height=2.5cm]{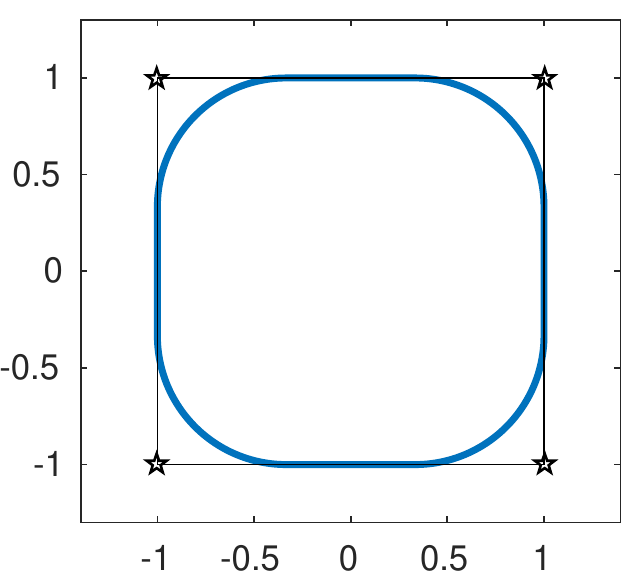} \hspace*{0.1cm}
   \includegraphics[height=2.5cm]{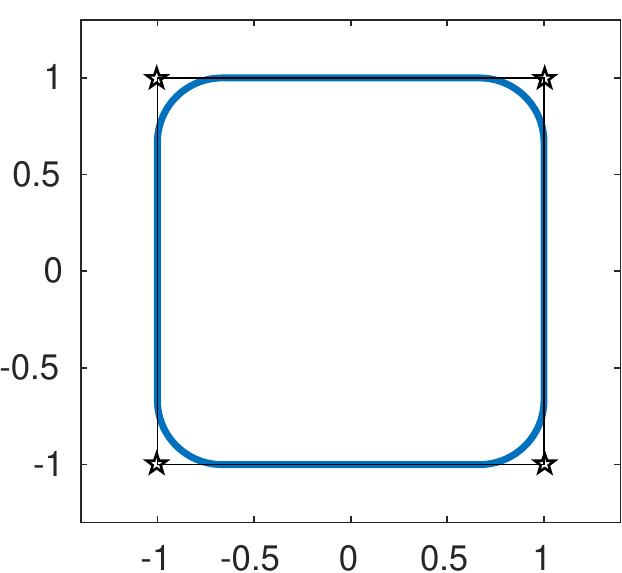}} 
  \caption{A $C^1$ smooth MDTB-spline representation of a square with  circular corners for $\ell=4^s$ and different values of $s$ according to Example~\ref{ex:example-curve}. The corresponding control points are indicated with the symbol $\star$.}\label{fig:example-curve}
\end{figure}

\begin{example}\label{ex:example-curve}
  Given $\ell\in(0,\infty)$, consider the periodic MDT-spline space $\splSpacerp$ defined by
  \begin{align*}
    \domain &= \left\{0, \frac{\pi}{2}, \ell+\frac{\pi}{2}, \ell+\pi, 2\ell+\pi,2\ell+\frac{3\pi}{2}, 3\ell+\frac{3\pi}{2},3\ell+2\pi,4\ell+2\pi\right\}, \\
    \bp &= \{2, 1, 2, 1, 2, 1, 2, 1\}, \quad
    \br = \{1,1,1,1,1,1,1,1,1\},
  \end{align*}
  and the local ECT-spaces
  \begin{equation*}
    \ECT{2}{(2i-1)} = \PT_2^{1} = \TT_2^{1}, \quad
    \ECT{1}{(2i)} = \PP_1, \quad
    i=1,2,3,4.
  \end{equation*}
  These local spaces can be represented as objects of the \Matlab{} classes \texttt{TB\_patch\_ptrig} (or, alternatively, \texttt{TB\_patch\_gtrig}) and \texttt{TB\_patch\_poly}.
  According to \eqref{eq:spline-dimension}, the dimension of the periodic MDT-spline is equal to $4$. The corresponding MDTB-spline basis functions are depicted in Figure~\ref{fig:example-curve}(a) for different values of $\ell=4^s$, $s=-2,-1,0,1$.  
  We now build a $C^1$ parametric spline curve whose parametric coefficients, called \emph{control points}, are given by
  \begin{equation*}
    (1,1), \quad (-1,1), \quad (-1,-1), \quad (1,-1).
  \end{equation*}
  The control polygon formed by these points is a square. The corresponding spline curve is depicted in Figure~\ref{fig:example-curve}(b) for different values of $\ell=4^s$, $s=-2,-1,0,1$. We clearly observe that it coincides with the control polygon except for the corners --- they are rounded and their size depends on $\ell$. In general, the explicit parametric expression of this curve is 
  \begin{equation*}
    (X(x),Y(x))=\begin{cases}
    \left(-L(2\sin(x)+\ell),\, L(2\cos(x)+\ell)\right), & x\in[0,\frac{\pi}{2}), \\
    \left(-1,\, -L(2x-\ell-\pi)\right), & x\in[\frac{\pi}{2}, \ell+\frac{\pi}{2}), \\
    \left(-L(2\sin(x-\ell)+\ell),\, L(2\cos(x-\ell)-\ell)\right), & x\in[\ell+\frac{\pi}{2}, \ell+\pi), \\
    \left(L(2x-3\ell-2\pi),\, -1\right), & x\in[\ell+\pi, 2\ell+\pi), \\
    \left(-L(2\sin(x-2\ell)-\ell),\, L(2\cos(x-2\ell)-\ell)\right), & x\in[2\ell+\pi, 2\ell+\frac{3\pi}{2}), \\
    \left(1,\, L(2x-5\ell-3\pi)\right), & x\in[2\ell+\frac{3\pi}{2}, 3\ell+\frac{3\pi}{2}), \\
    \left(-L(2\sin(x-3\ell)-\ell),\, L(2\cos(x-3\ell)+\ell)\right), & x\in[3\ell+\frac{3\pi}{2}, 3\ell+2\pi), \\
    \left(-L(2x-7\ell-4\pi),\, 1\right), & x\in[3\ell+2\pi, 4\ell+2\pi],
    \end{cases}
  \end{equation*}
  where $L:=1/(2+\ell)$. This expression shows that the radius of the four circular corners is equal to $2L$.
  In the limit cases, we obtain a perfect circle ($\ell=0$) or a perfect square ($\ell=\infty$).
  This example can be reproduced by executing the \Matlab{} script \texttt{EX\_basis\_B.m}. 
\end{example}

\subsection{Instabilities of ECT-Spaces} \label{sec:instability}

As mentioned in Section~\ref{sec:implementation-ECT}, the user of the \Matlab{} toolbox has to be aware of numerical instabilities in certain choices of ECT-spaces. We refer the reader to \citet{Roth:2019} for a detailed description and illustration of instabilities in the class of ECT-spaces described in Section~\ref{sec:ECT-Np}.

Possible sources of numerical instabilities are:
\begin{itemize}
  \item too small intervals $[\x_0,\x_1]$, so $\x_1\in(\x_0,\x_0+\varepsilon)$;
  \item roots close to each other, so $\croot_k=\alpha_k+\ii \beta_k$ and $\croot_{k+1}\in(\alpha_k-\varepsilon,\alpha_k+\varepsilon)+\ii (\beta_k-\varepsilon,\beta_k+\varepsilon)$;
  \item large exponential shape parameters, so $\alpha_k(\x_1-\x_0)\gg0$;
  \item high-dimensional spaces, so $\p\gg0$.
\end{itemize}
The first two bullets give rise to linear systems to be solved with similar (or almost identical) rows or columns in the conversion procedure described in Section~\ref{sec:computation-bernstein}, so these are (severely) ill-conditioned. The third bullet implies very high values of the function $\ee^{\alpha_k\x}$ and its derivatives. This leads to highly unbalanced entries in the involved linear systems, implying also ill-conditionings, and/or turns into numerical overflow.
The last bullet indicates that very high (end-point) derivatives are involved in the computation, and these will amplify the negative impact of the previous three sources.

Whenever a nearly singular linear system has to be solved in the \Matlab{} toolbox, a warning will be thrown. Then, it is up to the user to ignore the result or not as it may be inaccurate. 
The following checks are simple to perform and might give a numerical validation of the obtained result:
\begin{itemize}
  \item Check the non-negative partition of unity property of the computed \Chef{} Bernstein basis functions in a set of points distributed over the interval $[\x_0,\x_1]$. The basis functions can be simultaneously evaluated through the \Matlab{} method \texttt{TB\_evaluation\_all}. The resulting matrix should have non-negative entries and its column sum should give values (approximately) equal to one. 
  \item Do a visual inspection of the smoothness of the \Chef{} Bernstein basis functions; instabilities result in ``noisy'' functions. The basis functions can be simultaneously visualized through the \Matlab{} method \texttt{TB\_visualization\_all}.
\end{itemize}
As pointed out by \citet[Section~4.2]{Roth:2019}, given the infinitely many ECT-spaces with a vast possibility of inner structure, there is no general recipe  for the range of shape parameters for which the outputs of the proposed algorithms are guaranteed to be accurate enough --- this should be determined empirically by the user on a case-by-case basis. Nevertheless, a rule of thumb is to avoid high degrees and very small intervals.

\begin{figure}[t!]
  \subfigure[$\p=9$]
  {\includegraphics[height=4.4cm]{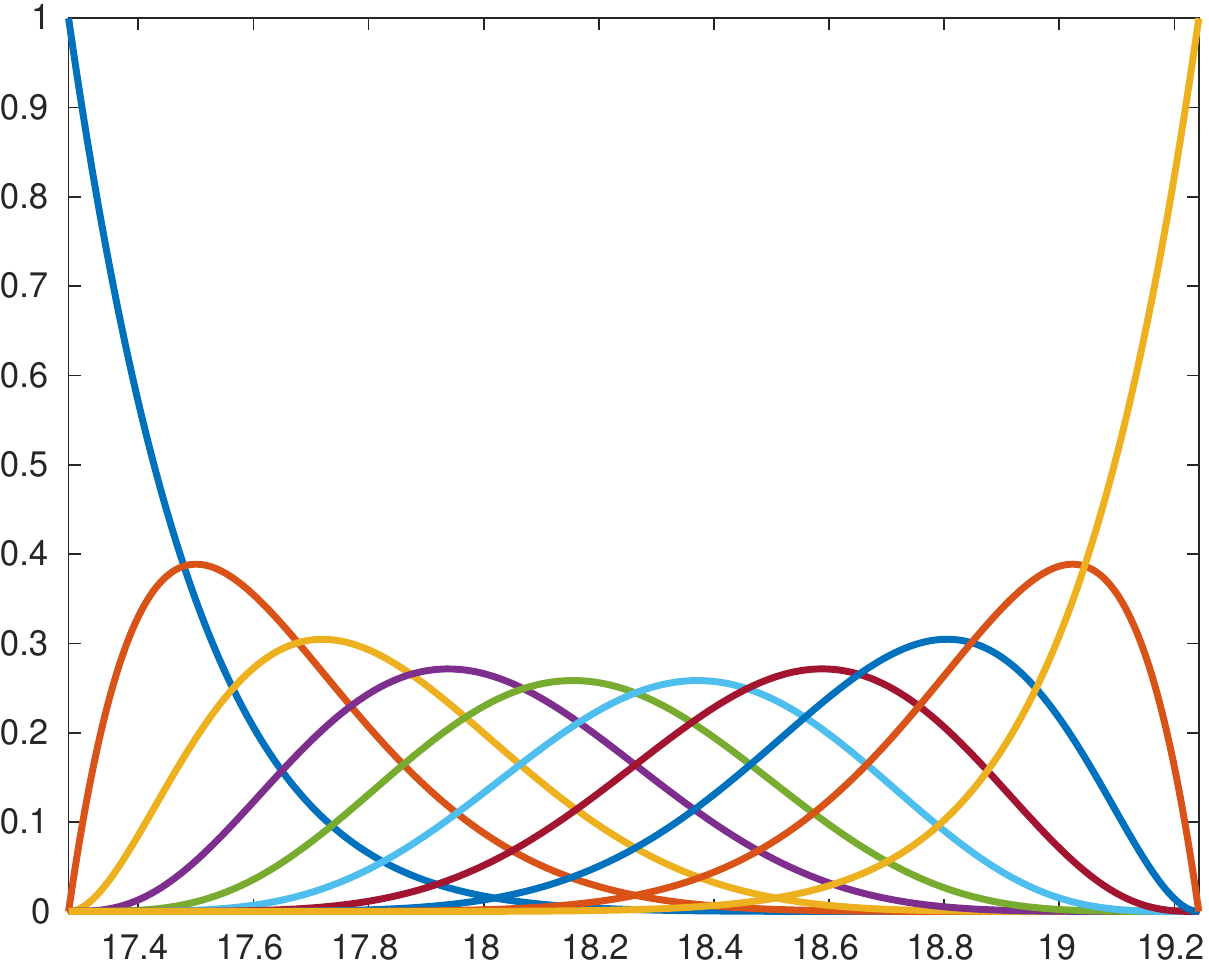}} \hspace*{0.5cm}
  \subfigure[$\p=10$]
  {\includegraphics[height=4.4cm]{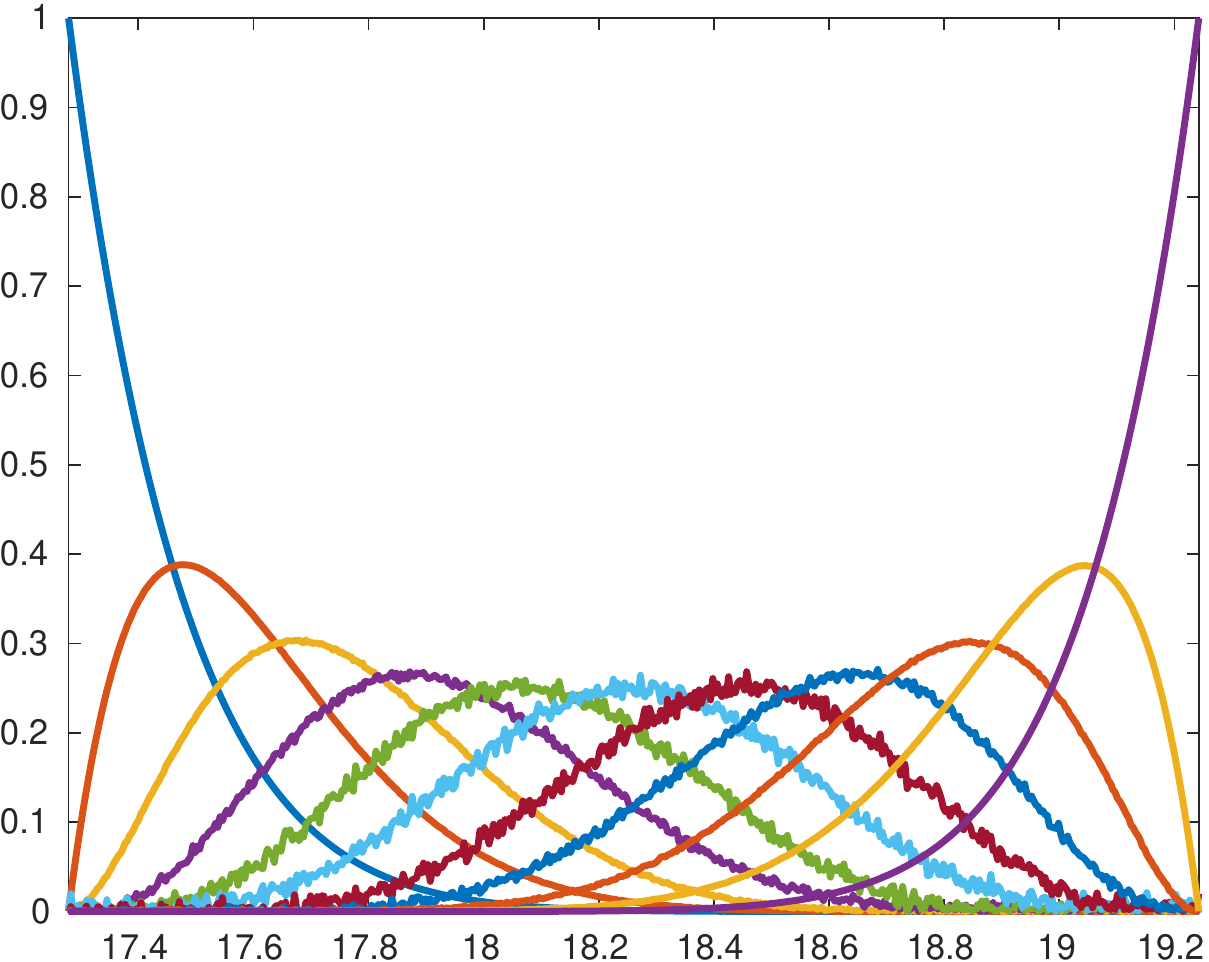}}
  \caption{Numerical instabilities in the computation of the \Chef{} Bernstein basis of the null-space $\WW_\p$ for high $\p$ according to Example~\ref{ex:instability-roth}.}\label{fig:instability-roth}
\end{figure}

\begin{example} \label{ex:instability-roth}
  Inspired by Example~2.7 from \citet{Roth:2019}, consider the null-space 
  \begin{equation*}
    \WW_\p:=\NN_\p^{(0,1,1),(\alpha_0,0,1),(\alpha_1,0,1),(\alpha_0,1,1)}, \quad \alpha_0:=\frac{1}{6\pi}, \quad \alpha_1:=\frac{1}{3\pi}
  \end{equation*}
  for $\p\geq 6$ on the interval $[\frac{11\pi}{2},\frac{49\pi}{8}]$. 
  The corresponding \Chef{} Bernstein basis functions can be computed using the general \Matlab{} class \texttt{TB\_patch\_tcheb}. They are visualized for $\p=9,10$ in Figure~\ref{fig:instability-roth}. In the case $\p=10$, one clearly notices that the central basis functions are perturbed by noise due to numerical instability.
  This is in agreement with the result obtained by \citet[Figure~9]{Roth:2019}. The partition-of-unity check, computed on a uniform grid of 501 points on $[\frac{11\pi}{2},\frac{49\pi}{8}]$, gives a maximal deviation from one 
  of $1.49\cdot10^{-4}$ for $\p=9$ and of $3.47\cdot10^{-2}$ for $\p=10$. 
  Higher values of $\p$ will lead to even more poor results.
  This example can be reproduced by executing the \Matlab{} script \texttt{EX\_instability\_A.m}.
\end{example}

For the generalized polynomial spaces $\EE_\p^\alpha$ and $\TT_\p^\beta$, a specialized implementation has been provided in the \Matlab{} classes \texttt{TB\_patch\_gexp} and \texttt{TB\_patch\_gtrig}, respectively. These implementations are more robust than the general null-space implementation of the class \texttt{TB\_patch\_tcheb}, and hence should be selected whenever possible. It is known that 
\begin{equation*}
  \lim_{\alpha\rightarrow0}\EE_\p^\alpha = \lim_{\beta\rightarrow0}\TT_\p^\beta =\PP_\p,
\end{equation*}
so for small values of the shape parameters the spaces $\EE_\p^\alpha$ and $\TT_\p^\beta$ behave like the algebraic polynomial space $\PP_\p$.
The same is true if, for fixed shape parameter, the length of the interval tends to zero. These cases can be a source of numerical instability in practice; they might arise in applications where mesh refinement and nested spaces are required, such as isogeometric analysis \cite{ManniRS:2017}. 
This problem has been properly treated in the specialized implementations as illustrated in the following example.

\begin{figure}[t!]
  \subfigure[\texttt{TB\_patch\_gtrig}]
  {\includegraphics[height=4.4cm]{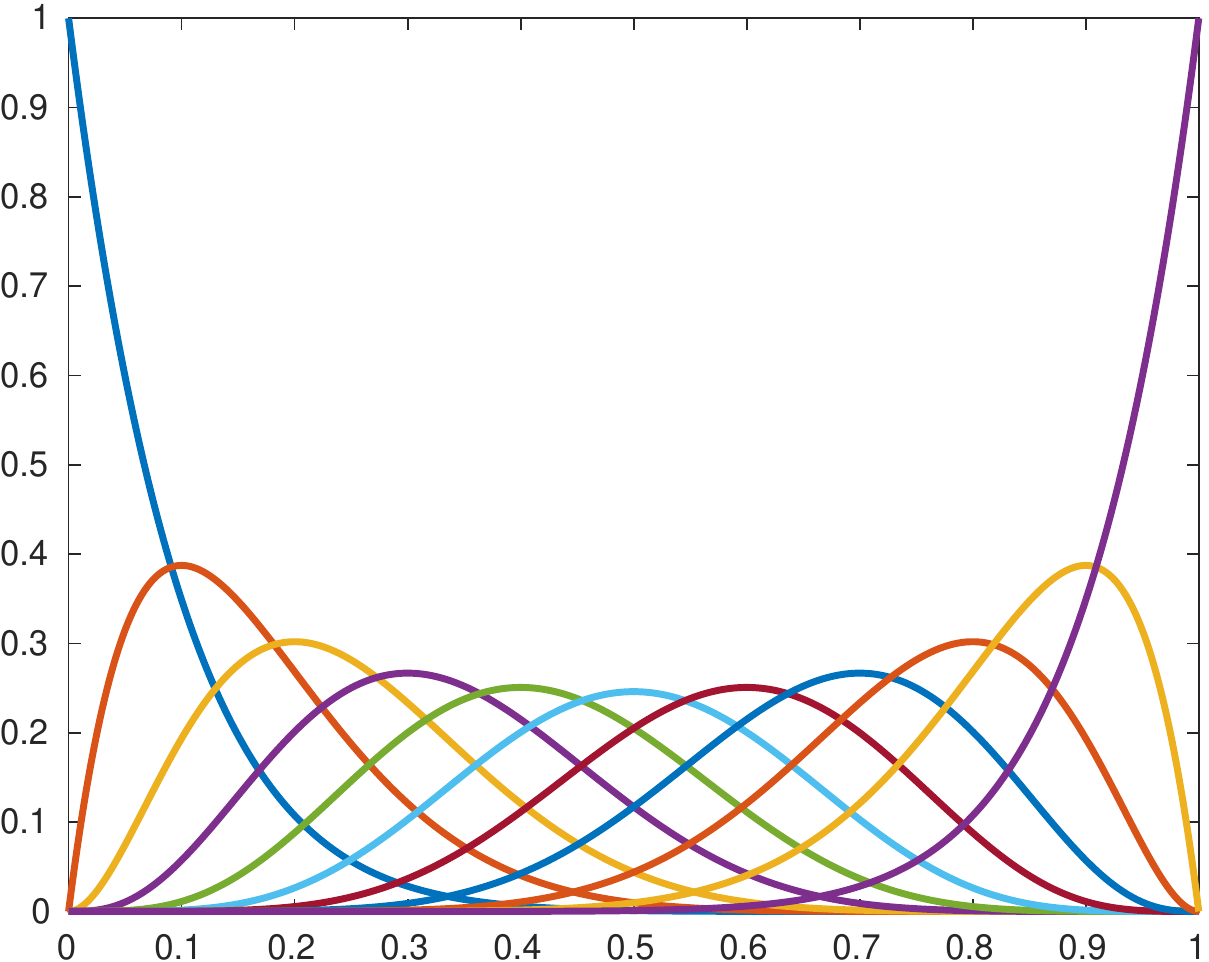}} \hspace*{0.5cm}
  \subfigure[\texttt{TB\_patch\_tcheb}]
  {\includegraphics[height=4.4cm]{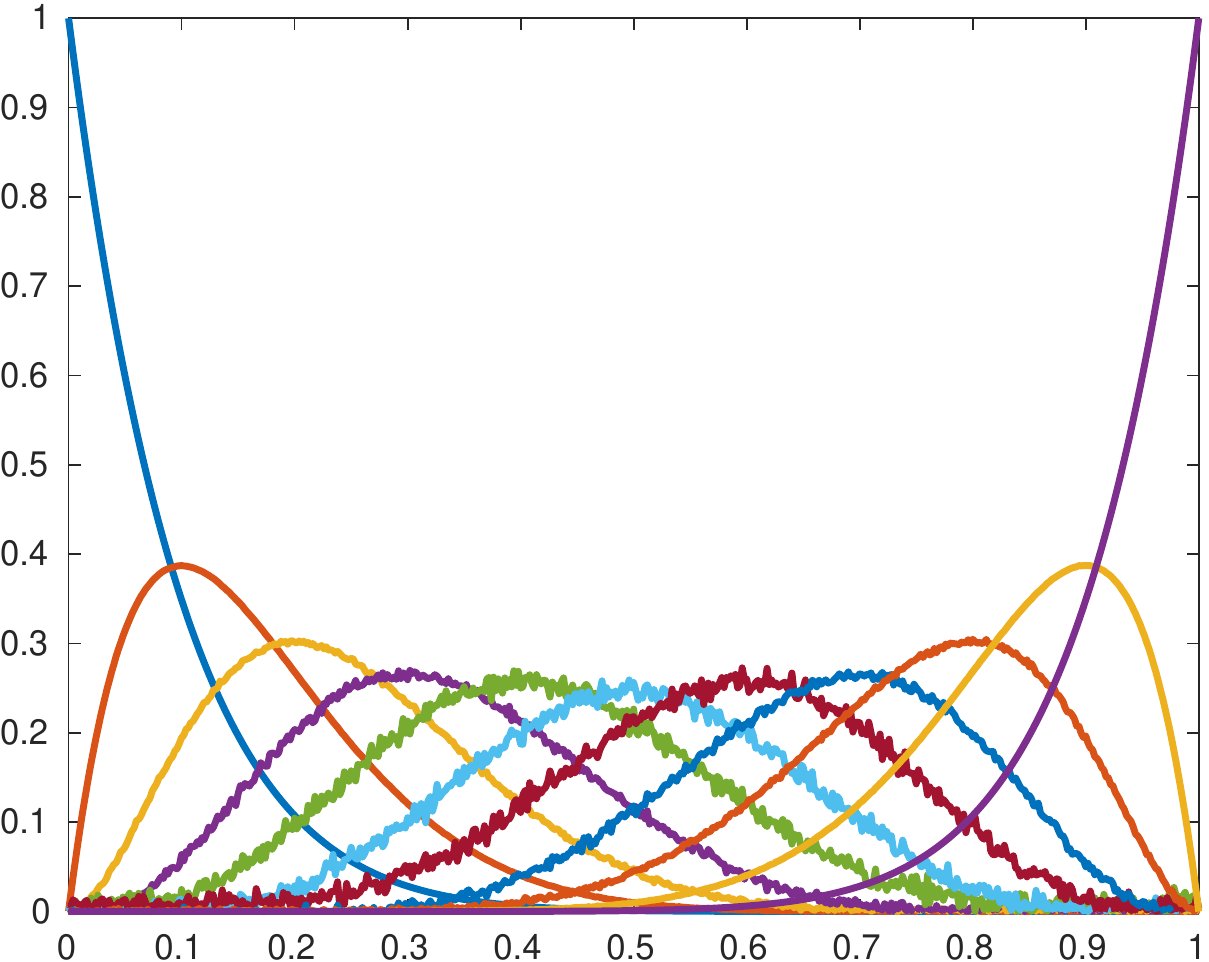}}
  \caption{The \Chef{} Bernstein basis of the generalized polynomial space $\TT_{10}^{1/3}$ computed by using different \Matlab{} classes in the toolbox according to Example~\ref{ex:instability-gtrig}.}\label{fig:instability-trig}
\end{figure}

\begin{example} \label{ex:instability-gtrig}
  Consider the generalized polynomial space $\TT_{10}^{1/3}$ on the interval $[0,1]$. We now compare the computation of the corresponding \Chef{} Bernstein basis functions using the \Matlab{} classes \texttt{TB\_patch\_gtrig} and \texttt{TB\_patch\_tcheb}. The outcome is visualized in Figure~\ref{fig:instability-trig}, and it is clear that the specialized class \texttt{TB\_patch\_gtrig} delivers a much more accurate result. This is confirmed by the partition-of-unity check; computed on a uniform grid of 501 points on $[0,1]$, it gives a maximal deviation from one 
  of $1.50\cdot10^{-10}$.
  This example can be reproduced by executing the \Matlab{} script \texttt{EX\_instability\_B.m}.
\end{example}

A comparable behavior is observed for the polynomial-type spaces $\PE_\p^\alpha$ and $\PT_\p^\beta$. The specialized \Matlab{} classes \texttt{TB\_patch\_pexp} and \texttt{TB\_patch\_ptrig}, respectively, provide a more robust implementation than the general null-space implementation of the class \texttt{TB\_patch\_tcheb}, and hence should also be selected whenever possible.

It is clear that working with MDT-spline spaces encounters similar problems of numerical instability as ECT-spaces, or even worse. The above proposed checks for a numerical validation can be applied in this context as well, using MDTB-splines instead of \Chef{} Bernstein functions. Moreover, one can add the following check:
\begin{itemize}
  \item Check that the extraction matrix has entries in the range $[0,1]$ and its column sum should give values (approximately) equal to one. This matrix can be computed through the \Matlab{} methods \texttt{MDTB\_extraction} or \texttt{MDTB\_extraction\_periodic}. 
\end{itemize}
As mentioned before, given the variety of possibilities of all parameters involved, it is impossible to provide a general recipe for their proper choice. Nevertheless, a rule of thumb is to use only stable local ECT-spaces and to avoid high degrees and highly non-uniform partitions.

\subsection{Critical Lengths for Design of ECT-Spaces} \label{sec:critical-length}

The choice of the interval $\interval:=[\x_0,\x_1]$ needs to be done in conjunction with the shape parameters in \eqref{eq:roots}, in order to ensure that the null-space $\NN_\p$ is an ECT-space possessing a \Chef{} Bernstein basis $\{B_{0,\p},\ldots,B_{\p,\p}\}$. 
The critical length for design $\ell'_\p$ guarantees the existence of such a basis, and so one has to choose $0<\x_1-\x_0<\ell'_\p$. The \Matlab{} toolbox works under the assumption that the interval of an ECT-space is always chosen within the critical length for design. 

In general, the determination of the critical length for design is theoretically challenging \cite{Brilleaud:2012,Carnicer:2003,CarnicerMP:2017} and one may need to rely on numerical techniques for its investigation \cite{BeccariCM:2020}. Our \Matlab{} toolbox could be of some help for this purpose as well. In the search for $\ell'_\p$, one could use the evaluation routines for the Bernstein basis to estimate the minimum positive value $\ell$ such that 
\begin{equation} \label{eq:crit-length-check}
\min_{0\leq j\leq \p}\min_{\x\in[0,\ell]}\bs_{j,\p}(\x)
\end{equation}
becomes negative. 
This condition can be approximately verified by tabulating the values of all Bernstein basis functions over a fine grid over $[0,\ell]$ for a sequence of increasing values of $\ell$. Once the minimum function value becomes negative (taking into account a negative tolerance), the previous value of $\ell$ in the sequence could be a possible candidate for the critical length for design. Of course, this approach assumes that there do not occur numerical instabilities in the computation of the Bernstein basis functions (see Section~\ref{sec:instability}). Moreover, 
it cannot give a conclusive result because the condition in \eqref{eq:crit-length-check} only provides an upper bound for $\ell'_\p$. Hence, the obtained numerical guess should still undergo a theoretical validation to be certain of the correct value of the critical length for design.

\begin{figure}[t!]
  \subfigure[Critical length $\ell'_\p$ of $\TT_\p^\beta$]
  {\includegraphics[height=4.4cm]{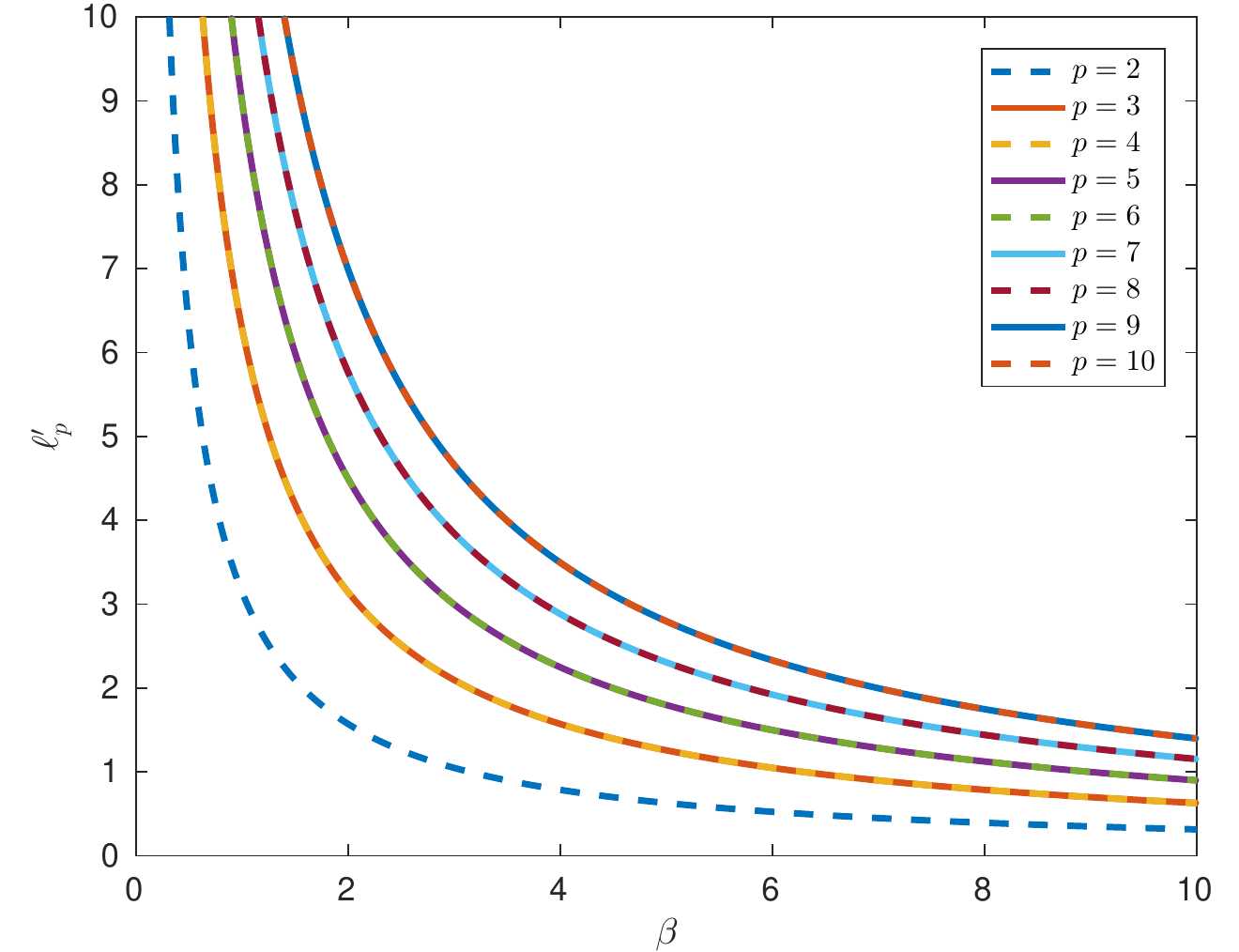}} \hspace*{0.3cm}
  \subfigure[Critical length $\ell'_6$ of $\NN_6^{\beta}$]
  {\includegraphics[height=4.4cm]{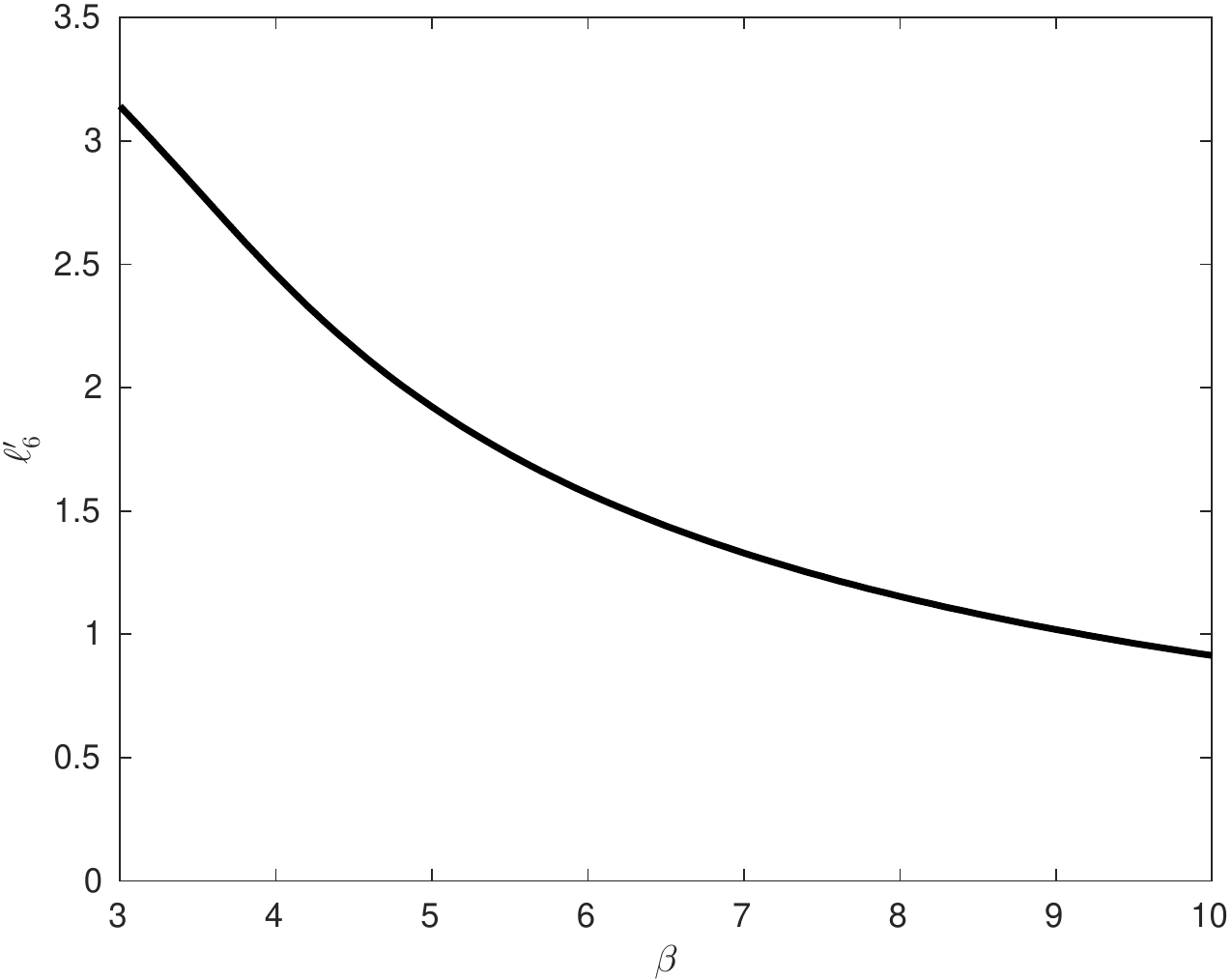}}
  \caption{Numerical computation of the critical length for design of the generalized polynomial space $\TT_\p^\beta$ for $\p\geq2$ and $\beta>0$ according to Example~\ref{ex:critical-length-trig} and
  of the null-space $\NN_6^{\beta}:=\NN_6^{(0,1,1),(0,2,1),(0,\beta,1)}$ for $\beta\geq3$ according to Example~\ref{ex:critical-length-mixed}.}\label{fig:critical-length-trig}
\end{figure}

\begin{example} \label{ex:critical-length-trig}
  Consider the generalized polynomial space $\TT_\p^\beta$ for $\p\geq2$ and $\beta=1$ on the interval $[0,\ell]$. The corresponding Bernstein basis functions can be evaluated using the \Matlab{} class \texttt{TB\_patch\_gtrig}. 
  For varying values of $\p$ and $\ell$, we compute the minimum value of the Bernstein basis functions, and arrive at the following numerical guesses of the critical lengths for design (up to three digits after the comma):
  \begin{equation} \label{eq:critical-length-trig}
    \tilde{\ell}'_2 = 3.141, \quad
    \tilde{\ell}'_3 = \tilde{\ell}'_4 = 6.283, \quad
    \tilde{\ell}'_5 = \tilde{\ell}'_6 = 8.986, \quad
    \tilde{\ell}'_7 = \tilde{\ell}'_8 = 11.526, \quad
    \tilde{\ell}'_9 = \tilde{\ell}'_{10} = 13.975.
  \end{equation}
  This is in agreement with the known values given in Example~\ref{ex:nullspace-trig}. Using a scaling argument, a numerical guess of the critical length for design of the space $\TT_\p^\beta$ for general $\beta>0$ can be obtained by dividing the values in \eqref{eq:critical-length-trig} by $\beta$.
  These are visualized as functions of $\beta\in(0,10]$ in Figure~\ref{fig:critical-length-trig}(a). This example can be reproduced by executing the \Matlab{} script \texttt{EX\_critical\_length\_A.m}.
\end{example}

\begin{example} \label{ex:critical-length-mixed}
  Consider the null-space $\NN_6^{\beta}:=\NN_6^{(0,1,1),(0,2,1),(0,\beta,1)}$ for $\beta\geq3$ on the interval $[0,\ell]$. 
  The corresponding Bernstein basis functions can be evaluated using the \Matlab{} class \texttt{TB\_patch\_tcheb}. 
  The minimum value of these functions is computed for varying $\beta$ and $\ell$, and our numerical guess of the critical length for design is visualized in Figure~\ref{fig:critical-length-trig}(b) as a function of $\beta\in[3,10]$.
  For $\beta=3$, we obtain the value $\tilde{\ell}'_6 = 3.141$ (up to three digits after the comma). In this case, we have $\NN_6^{3}=\PT_6^1$ and the critical length for design is known to be $\ell'_6=\pi$ \cite{Sanchez:1998}.
  This example can be reproduced by executing the \Matlab{} script \texttt{EX\_critical\_length\_B.m}.
\end{example}

Regarding MDT-spline spaces, the \Matlab{} toolbox works under the assumption that a valid MDTB-spline basis exists for the given user-specified parameters. We recall that the requirement of admissible weights in the sense of Definition~\ref{def:weight-assumption} is only a sufficient condition for the existence of MDTB-splines; see Remark~\ref{rmk:weight-assumption}.
We refer the reader to \citet[Remark~6.5]{Hiemstra:2020} for the specification of a more general setting. We also mention the work by \citet{Mazure:2011a} for a general characterization in terms of blossoms (in case of uniform degrees $\p_1=\cdots=\p_\nelms=\p$) and by \citet{BeccariCM:2019} for numerical tests. 
There is no direct extension of the concept of critical length for design towards MDT-spline spaces. However, one can possibly still provide a numerical check, similar to \eqref{eq:crit-length-check}, using MDTB-splines instead of \Chef{} Bernstein functions as illustrated in the following example. 

\begin{figure}[t!]
  \subfigure[Critical length $\ell'_6$ for $m=3$]
  {\includegraphics[height=4.4cm]{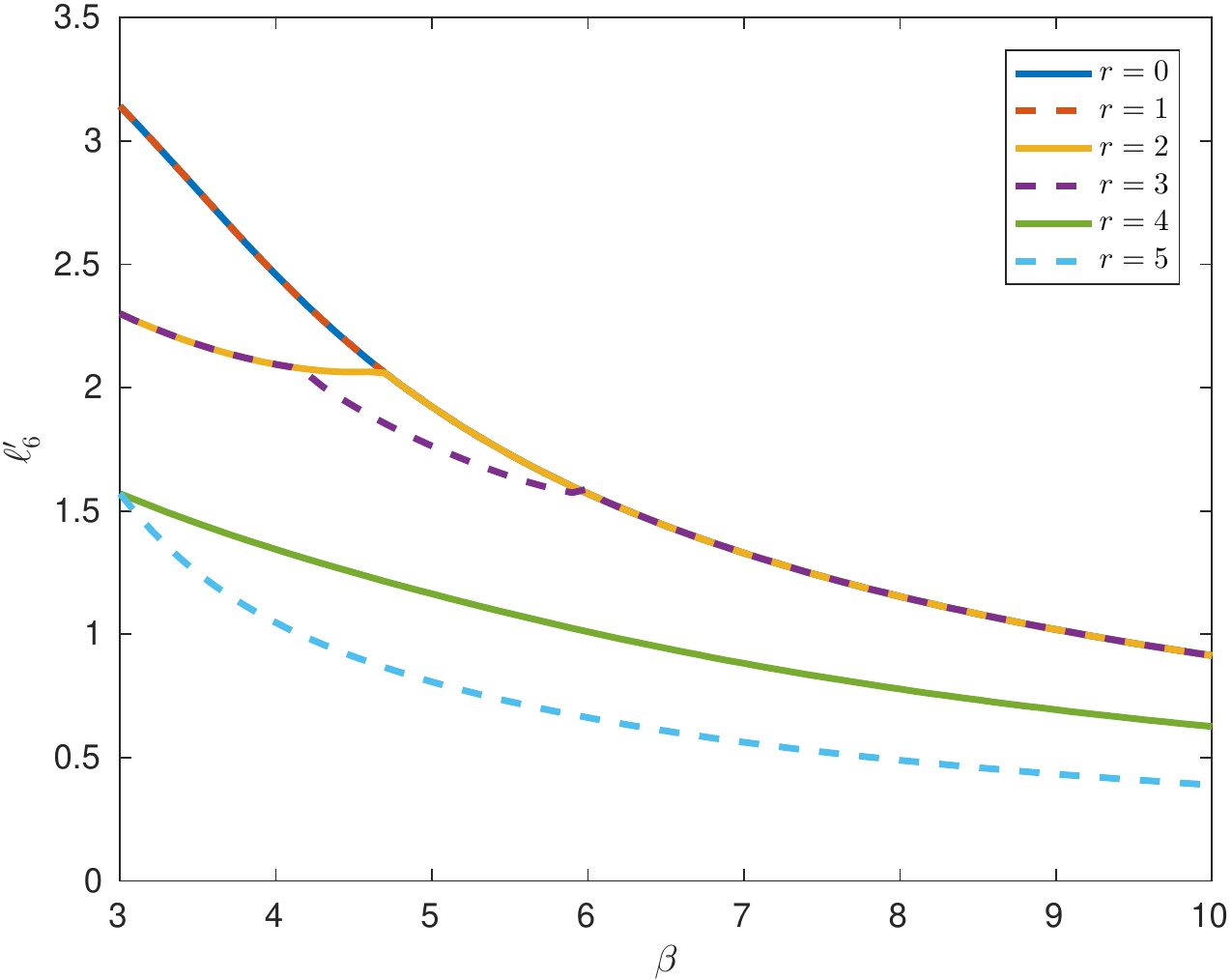}} \hspace*{0.3cm}
  \subfigure[Critical length $\ell'_6$ for $m=5$]
  {\includegraphics[height=4.4cm]{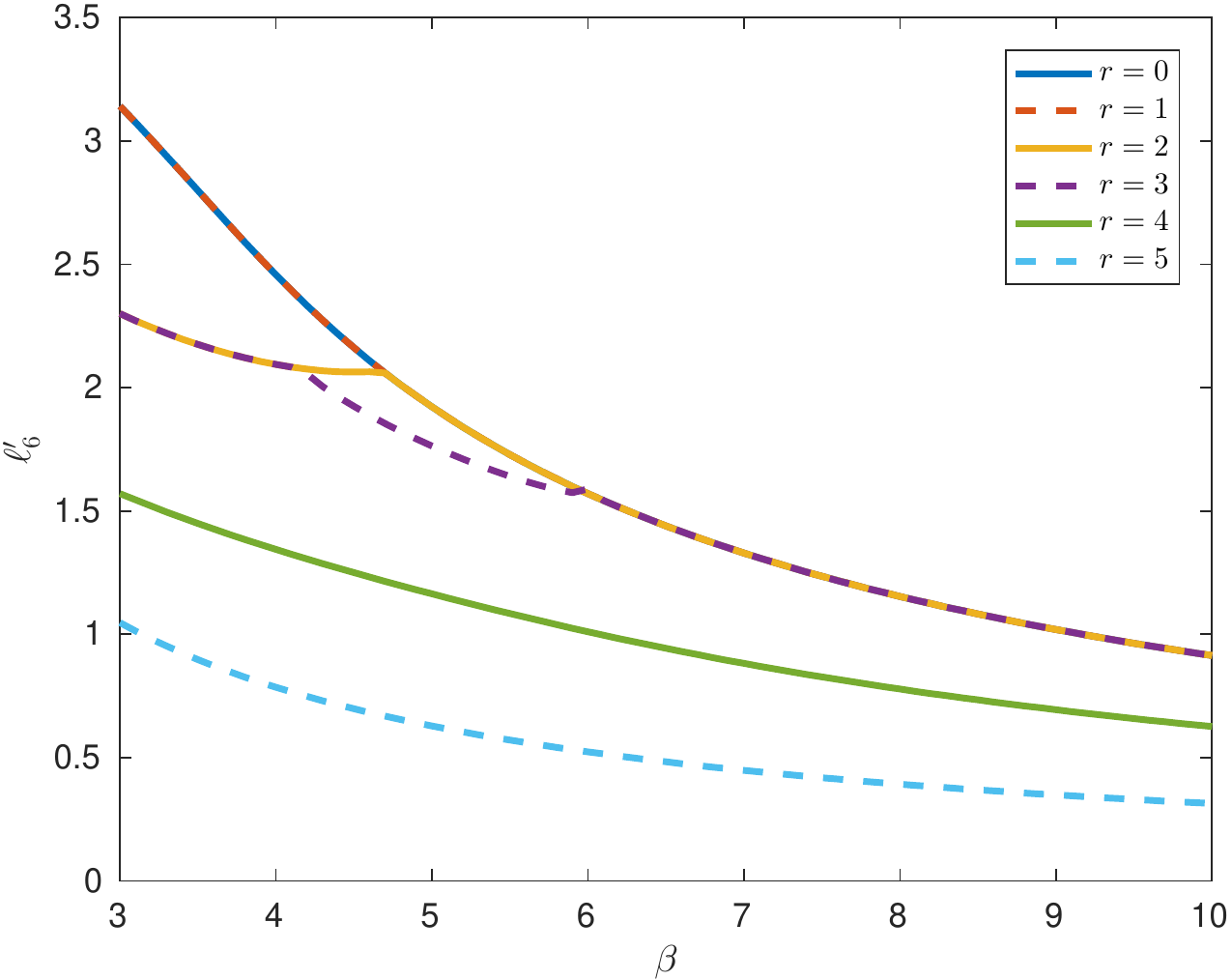}}
  \caption{Numerical computation of the critical length for design of the MDTB-spline space built locally of the null-space $\NN_6^{\beta}:=\NN_6^{(0,1,1),(0,2,1),(0,\beta,1)}$ for $\beta\geq3$ on $m$ subintervals of length $\ell$ with smoothness $\r\geq0$ according to Example~\ref{ex:critical-length-spline}.}\label{fig:critical-length-spline}
\end{figure}

\begin{example} \label{ex:critical-length-spline}
  As a continuation of Example~\ref{ex:critical-length-mixed}, consider now the uniform MDT-spline space $\splSpacerp$ built locally of the null-space $\NN_6^{\beta}:=\NN_6^{(0,1,1),(0,2,1),(0,\beta,1)}$ for $\beta\geq3$ on $m$ subintervals of length $\ell$, so 
  $\domain = \{0,\ell,2\ell,\ldots,m\ell\}$,
  and uniform smoothness $\br = \{-1,\r,\r,\ldots,\r,-1\}$ at the break points. For simplicity, also here we use the term \emph{critical length for design} to denote the supremum of the range of lengths $\ell>0$ that give rise to a valid MDT-spline space. The corresponding MDTB-splines can be obtained by using the factory class \texttt{MDTB\_patch\_tcheb}. 
  The minimum value of these functions is computed for $m=3,5$ and $\r=0,1,\ldots,5$, and varying $\beta$ and $\ell$. Our numerical guess of the critical length for design is visualized in Figure~\ref{fig:critical-length-spline} as a function of $\beta\in[3,10]$ for the different values of $m$ and $\r$. The result suggests that in the cases $\r=0,1$ the critical length for design is equal to the one of the local space $\NN_6^{\beta}$ (see Figure~\ref{fig:critical-length-trig}(b)), and in particular, equal to $\pi$ for $\beta=3$.
  This example can be reproduced by executing the \Matlab{} script \texttt{EX\_critical\_length\_C.m}.  
\end{example}

\section{Conclusion} \label{sec:conclusion}

In this article, we have presented a practical framework to deal with \Chef{} splines. These are splines with pieces drawn from different ECT-spaces (of possibly different dimensions) that are glued together smoothly. Under quite mild assumptions, they can be represented in terms of a so-called MDTB-spline basis, which enjoys properties similar to the classical polynomial B-spline basis. Thanks to the wide variety of ECT-spaces, such splines allow for an extraordinary flexibility in shape that may be optimally exploited in applications such as geometric modeling and isogeometric analysis.

We have detailed a simple procedure to built an extraction operator that represents all MDTB-splines as a linear combination of \Chef{} Bernstein functions related to each of the local ECT-spaces. With this procedure in hand, the complexity of computing and manipulating MDTB-splines is reduced to the same operations on \Chef{} Bernstein functions. The latter are, of course, easier to implement, but still require careful treatment in general as numerical instabilities are often lurking for certain (classes of) ECT-spaces; see Section~\ref{sec:instability}. 

Following this procedure, we have implemented an object-oriented \Matlab{} toolbox for dealing with MDTB-splines --- it is available through the CALGO library.
This is a redesigned and extended version of the MDB-spline toolbox by \citet{Speleers:2019}, developed for (algebraic) polynomial splines. The new toolbox can handle splines with pieces drawn from the large class of ECT-spaces that are null-spaces of constant-coefficient linear differential operators. This class covers most of the practically relevant ECT-spaces in geometric modeling \cite{Roth:2019} and isogeometric analysis \cite{ManniRS:2017}. Each of these spaces can be identified by means of the roots of the characteristic polynomial of the corresponding linear differential operator. Once all local ECT-spaces of the MDT-spline space of interest have been specified by the user, the toolbox allows for an easy computation and manipulation of the MDTB-splines through the \Matlab{} class \texttt{MDTB\_patch}.

At the same time, the treatment of ECT-spaces in the toolbox is valuable on its own, with the computation and manipulation of \Chef{} Bernstein functions. The general \Matlab{} class \texttt{TB\_patch\_tcheb} gives comparable results to the state-of-the-art C++ library recently developed by \citet{Roth:2019} --- the \Matlab{} implementation is slightly different however; see Remark~\ref{rmk:Roth}. In addition, the toolbox offers specialized (more efficient and/or more robust) implementations for algebraic polynomial spaces, other polynomial-type spaces, and generalized polynomial spaces --- important spaces in practical applications \cite{ManniRS:2017}. Thanks to the object-oriented structure of the toolbox, other implementations of ECT-spaces can be easily incorporated as well; one just needs to add new child classes of the abstract class \texttt{TB\_patch}. 

The \Matlab{} toolbox assumes that the input parameters specified by the user lead to valid ECT-spaces and MDT-spline spaces. For ECT-spaces, this means that the length of the basic interval should be smaller than its critical length for design. A sufficient condition for MDT-spline spaces is the existence of admissible weights in the sense of Definition~\ref{def:weight-assumption}. 
However, a full theoretical characterization of such spaces is not yet available in general. The development of practical recipes for selection of valid ranges of the input parameters is an important topic of further investigation. 
The toolbox could assist in this investigation; see Section~\ref{sec:critical-length}. 

The extraction procedure described in Section~\ref{sec:computation-spline}
is applicable to any kind of spline space that is equipped with a B-spline-like basis (in the sense of Remark~6.5 from \citet{Hiemstra:2020}), also beyond our \Chef{} setting. The \Matlab{} toolbox could thus be beneficial in this more general context as well. As long as child classes of the abstract class \texttt{TB\_patch} are provided, one can profit of the extraction procedure and the spline environment created by the toolbox.

\Chef{} splines are a beautiful theoretical tool, with a huge potential in applications such as geometric modeling and isogeometric analysis. However, as argued in \cite{Hiemstra:2020}, their full exploitation in practice has been limited so far by the lack of stable algorithms and implementations. We hope that the presented \Matlab{} toolbox may contribute towards the flourishing of \Chef{} splines, both in theory and practice.

\begin{acks}
This work was partially supported 
by the Beyond Borders Programme of the University of Rome Tor Vergata through the project ASTRID (CUP E84I19002250005) and 
by the MIUR Excellence Department Project awarded to the Department of Mathematics, University of Rome Tor Vergata (CUP E83C18000100006). 
The author is a member of Gruppo Nazionale per il Calcolo Scientifico --- Istituto Nazionale di Alta Matematica. 
\end{acks}


\bibliographystyle{ACM-Reference-Format}
\bibliography{mdtb_splines}


\begin{thebibliography}{44}


\ifx \showCODEN    \undefined \def \showCODEN     #1{\unskip}     \fi
\ifx \showDOI      \undefined \def \showDOI       #1{#1}\fi
\ifx \showISBNx    \undefined \def \showISBNx     #1{\unskip}     \fi
\ifx \showISBNxiii \undefined \def \showISBNxiii  #1{\unskip}     \fi
\ifx \showISSN     \undefined \def \showISSN      #1{\unskip}     \fi
\ifx \showLCCN     \undefined \def \showLCCN      #1{\unskip}     \fi
\ifx \shownote     \undefined \def \shownote      #1{#1}          \fi
\ifx \showarticletitle \undefined \def \showarticletitle #1{#1}   \fi
\ifx \showURL      \undefined \def \showURL       {\relax}        \fi
\providecommand\bibfield[2]{#2}
\providecommand\bibinfo[2]{#2}
\providecommand\natexlab[1]{#1}
\providecommand\showeprint[2][]{arXiv:#2}

\bibitem[\protect\citeauthoryear{Aimi, Diligenti, Sampoli, and Sestini}{Aimi
  et~al\mbox{.}}{2017}]%
        {Aimi:2017}
\bibfield{author}{\bibinfo{person}{A. Aimi}, \bibinfo{person}{M. Diligenti},
  \bibinfo{person}{M.~L. Sampoli}, {and} \bibinfo{person}{A. Sestini}.}
  \bibinfo{year}{2017}\natexlab{}.
\newblock \showarticletitle{Non-polynomial spline alternatives in isogeometric
  symmetric {G}alerkin {BEM}}.
\newblock \bibinfo{journal}{{\it Appl. Numer. Math.}}  \bibinfo{volume}{116}
  (\bibinfo{year}{2017}), \bibinfo{pages}{10--23}.
\newblock


\bibitem[\protect\citeauthoryear{Barry}{Barry}{1996}]%
        {Barry:1996}
\bibfield{author}{\bibinfo{person}{P.~J. Barry}.}
  \bibinfo{year}{1996}\natexlab{}.
\newblock \showarticletitle{de {B}oor--{F}ix dual functionals and algorithms
  for {T}chebycheffian {B}-spline curves}.
\newblock \bibinfo{journal}{{\it Constr. Approx.}}  \bibinfo{volume}{12}
  (\bibinfo{year}{1996}), \bibinfo{pages}{385--408}.
\newblock


\bibitem[\protect\citeauthoryear{Beccari and Casciola}{Beccari and
  Casciola}{2021}]%
        {Beccari:2021}
\bibfield{author}{\bibinfo{person}{C.~V. Beccari} {and} \bibinfo{person}{G.
  Casciola}.} \bibinfo{year}{2021}\natexlab{}.
\newblock \showarticletitle{Matrix representations for multi-degree
  {B}-splines}.
\newblock \bibinfo{journal}{{\it J. Comput. Appl. Math.}}
  \bibinfo{volume}{381}, Article \bibinfo{articleno}{113007}
  (\bibinfo{year}{2021}), \bibinfo{numpages}{18}~pages.
\newblock


\bibitem[\protect\citeauthoryear{Beccari, Casciola, and Mazure}{Beccari
  et~al\mbox{.}}{2019}]%
        {BeccariCM:2019}
\bibfield{author}{\bibinfo{person}{C.~V. Beccari}, \bibinfo{person}{G.
  Casciola}, {and} \bibinfo{person}{M.-L. Mazure}.}
  \bibinfo{year}{2019}\natexlab{}.
\newblock \showarticletitle{Design or not design? {A} numerical
  characterisation for piecewise {C}hebyshevian splines}.
\newblock \bibinfo{journal}{{\it Numer. Algorithms}}  \bibinfo{volume}{81}
  (\bibinfo{year}{2019}), \bibinfo{pages}{1--31}.
\newblock


\bibitem[\protect\citeauthoryear{Beccari, Casciola, and Mazure}{Beccari
  et~al\mbox{.}}{2020}]%
        {BeccariCM:2020}
\bibfield{author}{\bibinfo{person}{C.~V. Beccari}, \bibinfo{person}{G.
  Casciola}, {and} \bibinfo{person}{M.-L. Mazure}.}
  \bibinfo{year}{2020}\natexlab{}.
\newblock \showarticletitle{Critical length: An alternative approach}.
\newblock \bibinfo{journal}{{\it J. Comput. Appl. Math.}}
  \bibinfo{volume}{370}, Article \bibinfo{articleno}{112603}
  (\bibinfo{year}{2020}), \bibinfo{numpages}{16}~pages.
\newblock


\bibitem[\protect\citeauthoryear{Beccari, Casciola, and Morigi}{Beccari
  et~al\mbox{.}}{2017}]%
        {Beccari:2017}
\bibfield{author}{\bibinfo{person}{C.~V. Beccari}, \bibinfo{person}{G.
  Casciola}, {and} \bibinfo{person}{S. Morigi}.}
  \bibinfo{year}{2017}\natexlab{}.
\newblock \showarticletitle{On multi-degree splines}.
\newblock \bibinfo{journal}{{\it Comput. Aided Geom. Des.}}
  \bibinfo{volume}{58} (\bibinfo{year}{2017}), \bibinfo{pages}{8--23}.
\newblock


\bibitem[\protect\citeauthoryear{Bister and Prautzsch}{Bister and
  Prautzsch}{1997}]%
        {BisterP:1997}
\bibfield{author}{\bibinfo{person}{D. Bister} {and} \bibinfo{person}{H.
  Prautzsch}.} \bibinfo{year}{1997}\natexlab{}.
\newblock \showarticletitle{A new approach to {T}chebycheffian {B}-splines}.
\newblock In \bibinfo{booktitle}{{\em Curves and Surfaces with Applications in
  CAGD}}, \bibfield{editor}{\bibinfo{person}{A.~Le~M\'ehaut\'e},
  \bibinfo{person}{C.~Rabut}, {and} \bibinfo{person}{L.~L. Schumaker}} (Eds.).
  \bibinfo{publisher}{Vanderbilt University Press, Nashville},
  \bibinfo{pages}{387--394}.
\newblock


\bibitem[\protect\citeauthoryear{Brilleaud and Mazure}{Brilleaud and
  Mazure}{2012}]%
        {Brilleaud:2012}
\bibfield{author}{\bibinfo{person}{M. Brilleaud} {and} \bibinfo{person}{M.-L.
  Mazure}.} \bibinfo{year}{2012}\natexlab{}.
\newblock \showarticletitle{Mixed hyperbolic/trigonometric spaces for design}.
\newblock \bibinfo{journal}{{\it Comput. Math. Appl.}}  \bibinfo{volume}{64}
  (\bibinfo{year}{2012}), \bibinfo{pages}{2459--2477}.
\newblock


\bibitem[\protect\citeauthoryear{Buchwald and M\"uhlbach}{Buchwald and
  M\"uhlbach}{2003}]%
        {Buchwald:2003}
\bibfield{author}{\bibinfo{person}{B. Buchwald} {and} \bibinfo{person}{G.
  M\"uhlbach}.} \bibinfo{year}{2003}\natexlab{}.
\newblock \showarticletitle{Construction of {B}-splines for generalized spline
  spaces generated from local {ECT}-systems}.
\newblock \bibinfo{journal}{{\it J. Comput. Appl. Math.}}
  \bibinfo{volume}{159} (\bibinfo{year}{2003}), \bibinfo{pages}{249--267}.
\newblock


\bibitem[\protect\citeauthoryear{Carnicer, Mainar, and Pe\~na}{Carnicer
  et~al\mbox{.}}{2003}]%
        {Carnicer:2003}
\bibfield{author}{\bibinfo{person}{J.~M. Carnicer}, \bibinfo{person}{E.
  Mainar}, {and} \bibinfo{person}{J.~M. Pe\~na}.}
  \bibinfo{year}{2003}\natexlab{}.
\newblock \showarticletitle{Critical length for design purposes and extended
  {C}hebyshev spaces}.
\newblock \bibinfo{journal}{{\it Constr. Approx.}}  \bibinfo{volume}{20}
  (\bibinfo{year}{2003}), \bibinfo{pages}{55--71}.
\newblock


\bibitem[\protect\citeauthoryear{Carnicer, Mainar, and Pe\~na}{Carnicer
  et~al\mbox{.}}{2017}]%
        {CarnicerMP:2017}
\bibfield{author}{\bibinfo{person}{J.~M. Carnicer}, \bibinfo{person}{E.
  Mainar}, {and} \bibinfo{person}{J.~M. Pe\~na}.}
  \bibinfo{year}{2017}\natexlab{}.
\newblock \showarticletitle{Critical lengths of cycloidal spaces are zeros of
  {B}essel functions}.
\newblock \bibinfo{journal}{{\em Calcolo\/}}  \bibinfo{volume}{54}
  (\bibinfo{year}{2017}), \bibinfo{pages}{1521--1531}.
\newblock


\bibitem[\protect\citeauthoryear{Cohen, Riesenfeld, and Elber}{Cohen
  et~al\mbox{.}}{2001}]%
        {Cohen:2001}
\bibfield{author}{\bibinfo{person}{E. Cohen}, \bibinfo{person}{R.~F.
  Riesenfeld}, {and} \bibinfo{person}{G. Elber}.}
  \bibinfo{year}{2001}\natexlab{}.
\newblock \bibinfo{booktitle}{{\em Geometric Modeling with Splines: An
  Introduction}}.
\newblock \bibinfo{publisher}{CRC Press}.
\newblock


\bibitem[\protect\citeauthoryear{Coppel}{Coppel}{1971}]%
        {Coppel:1971}
\bibfield{author}{\bibinfo{person}{W.~A. Coppel}.}
  \bibinfo{year}{1971}\natexlab{}.
\newblock \bibinfo{booktitle}{{\em Disconjugacy}}.
\newblock \bibinfo{publisher}{Springer-Verlag}.
\newblock


\bibitem[\protect\citeauthoryear{Costantini, Lyche, and Manni}{Costantini
  et~al\mbox{.}}{2005}]%
        {Costantini:2005}
\bibfield{author}{\bibinfo{person}{P. Costantini}, \bibinfo{person}{T. Lyche},
  {and} \bibinfo{person}{C. Manni}.} \bibinfo{year}{2005}\natexlab{}.
\newblock \showarticletitle{On a class of weak {T}chebycheff systems}.
\newblock \bibinfo{journal}{{\it Numer. Math.}}  \bibinfo{volume}{101}
  (\bibinfo{year}{2005}), \bibinfo{pages}{333--354}.
\newblock


\bibitem[\protect\citeauthoryear{Cottrell, Hughes, and Bazilevs}{Cottrell
  et~al\mbox{.}}{2009}]%
        {Cottrell:2009}
\bibfield{author}{\bibinfo{person}{J.~A. Cottrell}, \bibinfo{person}{T.~J.~R.
  Hughes}, {and} \bibinfo{person}{Y. Bazilevs}.}
  \bibinfo{year}{2009}\natexlab{}.
\newblock \bibinfo{booktitle}{{\em Isogeometric Analysis: Toward Integration of
  {CAD} and {FEA}}}.
\newblock \bibinfo{publisher}{John Wiley \& Sons}.
\newblock


\bibitem[\protect\citeauthoryear{Fang, Ma, and Wang}{Fang
  et~al\mbox{.}}{2010}]%
        {Fang:2010}
\bibfield{author}{\bibinfo{person}{M. Fang}, \bibinfo{person}{W. Ma}, {and}
  \bibinfo{person}{G. Wang}.} \bibinfo{year}{2010}\natexlab{}.
\newblock \showarticletitle{A generalized curve subdivision scheme of arbitrary
  order with a tension parameter}.
\newblock \bibinfo{journal}{{\it Comput. Aided Geom. Des.}}
  \bibinfo{volume}{27} (\bibinfo{year}{2010}), \bibinfo{pages}{720--733}.
\newblock


\bibitem[\protect\citeauthoryear{Hiemstra, R., Manni, Speleers, and
  Toshniwal}{Hiemstra et~al\mbox{.}}{2020}]%
        {Hiemstra:2020}
\bibfield{author}{\bibinfo{person}{R.~R. Hiemstra}, \bibinfo{person}{Hughes
  T.~J. R.}, \bibinfo{person}{C. Manni}, \bibinfo{person}{H. Speleers}, {and}
  \bibinfo{person}{D. Toshniwal}.} \bibinfo{year}{2020}\natexlab{}.
\newblock \showarticletitle{A {T}chebycheffian extension of multi-degree
  {B}-splines: Algorithmic computation and properties}.
\newblock \bibinfo{journal}{{\it SIAM J. Numer. Anal.}}  \bibinfo{volume}{58}
  (\bibinfo{year}{2020}), \bibinfo{pages}{1138--1163}.
\newblock


\bibitem[\protect\citeauthoryear{Karlin}{Karlin}{1968}]%
        {Karlin:1968}
\bibfield{author}{\bibinfo{person}{S. Karlin}.}
  \bibinfo{year}{1968}\natexlab{}.
\newblock \bibinfo{booktitle}{{\em Total Positivity}}.
\newblock \bibinfo{publisher}{Stanford University Press}.
\newblock


\bibitem[\protect\citeauthoryear{Karlin and Ziegler}{Karlin and
  Ziegler}{1966}]%
        {KarlinZ:1966}
\bibfield{author}{\bibinfo{person}{S. Karlin} {and} \bibinfo{person}{Z.
  Ziegler}.} \bibinfo{year}{1966}\natexlab{}.
\newblock \showarticletitle{{C}hebyshevian spline functions}.
\newblock \bibinfo{journal}{{\it SIAM J. Numer. Anal.}}  \bibinfo{volume}{3}
  (\bibinfo{year}{1966}), \bibinfo{pages}{514--543}.
\newblock


\bibitem[\protect\citeauthoryear{Koch and Lyche}{Koch and Lyche}{1993}]%
        {KochL1993}
\bibfield{author}{\bibinfo{person}{P.~E. Koch} {and} \bibinfo{person}{T.
  Lyche}.} \bibinfo{year}{1993}\natexlab{}.
\newblock \showarticletitle{Interpolation with exponential {B}-splines in
  tension}.
\newblock In \bibinfo{booktitle}{{\em Geometric Modelling}},
  \bibfield{editor}{\bibinfo{person}{G.~Farin}, \bibinfo{person}{H.~Hagen},
  \bibinfo{person}{H.~Noltemeier}, {and} \bibinfo{person}{W.~Kn\"{o}del}}
  (Eds.). \bibinfo{publisher}{Springer--Verlag, Wien},
  \bibinfo{pages}{173--190}.
\newblock


\bibitem[\protect\citeauthoryear{Kvasov and Sattayatham}{Kvasov and
  Sattayatham}{1999}]%
        {KvasovS:1999}
\bibfield{author}{\bibinfo{person}{B. Kvasov} {and} \bibinfo{person}{P.
  Sattayatham}.} \bibinfo{year}{1999}\natexlab{}.
\newblock \showarticletitle{{GB}-splines of arbitrary order}.
\newblock \bibinfo{journal}{{\it J. Comput. Appl. Math.}}
  \bibinfo{volume}{104} (\bibinfo{year}{1999}), \bibinfo{pages}{63--88}.
\newblock


\bibitem[\protect\citeauthoryear{Lyche}{Lyche}{1985}]%
        {Lyche:1985}
\bibfield{author}{\bibinfo{person}{T. Lyche}.} \bibinfo{year}{1985}\natexlab{}.
\newblock \showarticletitle{A recurrence relation for {C}hebyshevian
  {B}-splines}.
\newblock \bibinfo{journal}{{\it Constr. Approx.}}  \bibinfo{volume}{1}
  (\bibinfo{year}{1985}), \bibinfo{pages}{155--173}.
\newblock


\bibitem[\protect\citeauthoryear{Lyche, Manni, and Speleers}{Lyche
  et~al\mbox{.}}{2019}]%
        {Lyche:2019}
\bibfield{author}{\bibinfo{person}{T. Lyche}, \bibinfo{person}{C. Manni}, {and}
  \bibinfo{person}{H. Speleers}.} \bibinfo{year}{2019}\natexlab{}.
\newblock \showarticletitle{Tchebycheffian {B}-splines revisited: An
  introductory exposition}.
\newblock In \bibinfo{booktitle}{{\em Advanced Methods for Geometric Modeling
  and Numerical Simulation}}, \bibfield{editor}{\bibinfo{person}{C.~Giannelli}
  {and} \bibinfo{person}{H.~Speleers}} (Eds.). \bibinfo{series}{Springer INdAM
  Series}, Vol.~\bibinfo{volume}{35}. \bibinfo{publisher}{Springer
  International Publishing AG}, \bibinfo{pages}{179--216}.
\newblock


\bibitem[\protect\citeauthoryear{Manni, Pelosi, and Sampoli}{Manni
  et~al\mbox{.}}{2011}]%
        {ManniPS:2011}
\bibfield{author}{\bibinfo{person}{C. Manni}, \bibinfo{person}{F. Pelosi},
  {and} \bibinfo{person}{M.~L. Sampoli}.} \bibinfo{year}{2011}\natexlab{}.
\newblock \showarticletitle{Generalized {B}-splines as a tool in isogeometric
  analysis}.
\newblock \bibinfo{journal}{{\it Comput. Methods Appl. Mech. Eng.}}
  \bibinfo{volume}{200} (\bibinfo{year}{2011}), \bibinfo{pages}{867--881}.
\newblock


\bibitem[\protect\citeauthoryear{Manni, Reali, and Speleers}{Manni
  et~al\mbox{.}}{2015}]%
        {ManniRS:2015}
\bibfield{author}{\bibinfo{person}{C. Manni}, \bibinfo{person}{A. Reali}, {and}
  \bibinfo{person}{H. Speleers}.} \bibinfo{year}{2015}\natexlab{}.
\newblock \showarticletitle{Isogeometric collocation methods with generalized
  B-splines}.
\newblock \bibinfo{journal}{{\it Comput. Math. Appl.}}  \bibinfo{volume}{70}
  (\bibinfo{year}{2015}), \bibinfo{pages}{1659--1675}.
\newblock


\bibitem[\protect\citeauthoryear{Manni, Roman, and Speleers}{Manni
  et~al\mbox{.}}{2017}]%
        {ManniRS:2017}
\bibfield{author}{\bibinfo{person}{C. Manni}, \bibinfo{person}{F. Roman}, {and}
  \bibinfo{person}{H. Speleers}.} \bibinfo{year}{2017}\natexlab{}.
\newblock \showarticletitle{Generalized {B}-splines in isogeometric analysis}.
\newblock In \bibinfo{booktitle}{{\em Approximation Theory XV: San Antonio
  2016}}, \bibfield{editor}{\bibinfo{person}{G.~E. Fasshauer} {and}
  \bibinfo{person}{L.~L. Schumaker}} (Eds.). \bibinfo{series}{Springer
  Proceedings in Mathematics \& Statistics}, Vol.~\bibinfo{volume}{201}.
  \bibinfo{publisher}{Springer International Publishing AG},
  \bibinfo{pages}{239--267}.
\newblock


\bibitem[\protect\citeauthoryear{Mazure}{Mazure}{2007}]%
        {Mazure:2007}
\bibfield{author}{\bibinfo{person}{M.-L. Mazure}.}
  \bibinfo{year}{2007}\natexlab{}.
\newblock \showarticletitle{Extended {C}hebyshev piecewise spaces characterised
  via weight functions}.
\newblock \bibinfo{journal}{{\it J. Approx. Theory}}  \bibinfo{volume}{145}
  (\bibinfo{year}{2007}), \bibinfo{pages}{33--54}.
\newblock


\bibitem[\protect\citeauthoryear{Mazure}{Mazure}{2011a}]%
        {Mazure:2011b}
\bibfield{author}{\bibinfo{person}{M.-L. Mazure}.}
  \bibinfo{year}{2011}\natexlab{a}.
\newblock \showarticletitle{Finding all systems of weight functions associated
  with a given extended {C}hebyshev space}.
\newblock \bibinfo{journal}{{\it J. Approx. Theory}}  \bibinfo{volume}{163}
  (\bibinfo{year}{2011}), \bibinfo{pages}{363--376}.
\newblock


\bibitem[\protect\citeauthoryear{Mazure}{Mazure}{2011b}]%
        {Mazure:2011a}
\bibfield{author}{\bibinfo{person}{M.-L. Mazure}.}
  \bibinfo{year}{2011}\natexlab{b}.
\newblock \showarticletitle{How to build all {C}hebyshevian spline spaces good
  for geometric design?}
\newblock \bibinfo{journal}{{\it Numer. Math.}}  \bibinfo{volume}{119}
  (\bibinfo{year}{2011}), \bibinfo{pages}{517--556}.
\newblock


\bibitem[\protect\citeauthoryear{Mazure}{Mazure}{2018}]%
        {Mazure:2018}
\bibfield{author}{\bibinfo{person}{M.-L. Mazure}.}
  \bibinfo{year}{2018}\natexlab{}.
\newblock \showarticletitle{Constructing totally positive piecewise
  {C}hebyshevian {B}-spline bases}.
\newblock \bibinfo{journal}{{\it J. Comput. Appl. Math.}}
  \bibinfo{volume}{342} (\bibinfo{year}{2018}), \bibinfo{pages}{550--586}.
\newblock


\bibitem[\protect\citeauthoryear{N\"urnberger, Schumaker, Sommer, and
  Strauss}{N\"urnberger et~al\mbox{.}}{1983}]%
        {Nurnberger:1983}
\bibfield{author}{\bibinfo{person}{G. N\"urnberger}, \bibinfo{person}{L.~L.
  Schumaker}, \bibinfo{person}{M. Sommer}, {and} \bibinfo{person}{H. Strauss}.}
  \bibinfo{year}{1983}\natexlab{}.
\newblock \showarticletitle{Interpolation by generalized splines}.
\newblock \bibinfo{journal}{{\it Numer. Math.}}  \bibinfo{volume}{42}
  (\bibinfo{year}{1983}), \bibinfo{pages}{195--212}.
\newblock


\bibitem[\protect\citeauthoryear{N\"urnberger, Schumaker, Sommer, and
  Strauss}{N\"urnberger et~al\mbox{.}}{1984}]%
        {Nurnberger:1984}
\bibfield{author}{\bibinfo{person}{G. N\"urnberger}, \bibinfo{person}{L.~L.
  Schumaker}, \bibinfo{person}{M. Sommer}, {and} \bibinfo{person}{H. Strauss}.}
  \bibinfo{year}{1984}\natexlab{}.
\newblock \showarticletitle{Generalized {C}hebyshevian splines}.
\newblock \bibinfo{journal}{{\it SIAM J. Math. Anal.}}  \bibinfo{volume}{15}
  (\bibinfo{year}{1984}), \bibinfo{pages}{790--804}.
\newblock


\bibitem[\protect\citeauthoryear{Pottmann}{Pottmann}{1993}]%
        {Pottmann:1993}
\bibfield{author}{\bibinfo{person}{H. Pottmann}.}
  \bibinfo{year}{1993}\natexlab{}.
\newblock \showarticletitle{The geometry of {T}chebycheffian splines}.
\newblock \bibinfo{journal}{{\it Comput. Aided Geom. Des.}}
  \bibinfo{volume}{10} (\bibinfo{year}{1993}), \bibinfo{pages}{181--210}.
\newblock


\bibitem[\protect\citeauthoryear{Roman, Manni, and Speleers}{Roman
  et~al\mbox{.}}{2017}]%
        {RomanMS:2017}
\bibfield{author}{\bibinfo{person}{F. Roman}, \bibinfo{person}{C. Manni}, {and}
  \bibinfo{person}{H. Speleers}.} \bibinfo{year}{2017}\natexlab{}.
\newblock \showarticletitle{Numerical approximation of {GB}-splines by a
  convolutional approach}.
\newblock \bibinfo{journal}{{\it Appl. Numer. Math.}}  \bibinfo{volume}{116}
  (\bibinfo{year}{2017}), \bibinfo{pages}{273--285}.
\newblock


\bibitem[\protect\citeauthoryear{R\'oth}{R\'oth}{2019}]%
        {Roth:2019}
\bibfield{author}{\bibinfo{person}{\'A. R\'oth}.}
  \bibinfo{year}{2019}\natexlab{}.
\newblock \showarticletitle{Algorithm 992: An {OpenGL}- and {C++}-based
  function library for curve and surface modeling in a large class of extended
  {C}hebyshev spaces}.
\newblock \bibinfo{journal}{{\it ACM Trans. Math. Software}}
  \bibinfo{volume}{45}, Article \bibinfo{articleno}{13} (\bibinfo{year}{2019}),
  \bibinfo{numpages}{32}~pages.
\newblock


\bibitem[\protect\citeauthoryear{S\'anchez-Reyes}{S\'anchez-Reyes}{1998}]%
        {Sanchez:1998}
\bibfield{author}{\bibinfo{person}{J. S\'anchez-Reyes}.}
  \bibinfo{year}{1998}\natexlab{}.
\newblock \showarticletitle{Harmonic rational {B}\'ezier curves, p-{B}\'ezier
  curves and trigonometric polynomials}.
\newblock \bibinfo{journal}{{\it Comput. Aided Geom. Des.}}
  \bibinfo{volume}{15} (\bibinfo{year}{1998}), \bibinfo{pages}{909--923}.
\newblock


\bibitem[\protect\citeauthoryear{Schumaker}{Schumaker}{2007}]%
        {Schumaker:2007}
\bibfield{author}{\bibinfo{person}{L.~L. Schumaker}.}
  \bibinfo{year}{2007}\natexlab{}.
\newblock \bibinfo{booktitle}{{\em Spline Functions: Basic Theory, Third
  Edition}}.
\newblock \bibinfo{publisher}{Cambridge University Press}.
\newblock


\bibitem[\protect\citeauthoryear{Shen and Wang}{Shen and Wang}{2005}]%
        {Shen:2005}
\bibfield{author}{\bibinfo{person}{W.-Q. Shen} {and} \bibinfo{person}{G.-Z.
  Wang}.} \bibinfo{year}{2005}\natexlab{}.
\newblock \showarticletitle{A class of quasi B\'ezier curves based on
  hyperbolic polynomials}.
\newblock \bibinfo{journal}{{\it J. Zhejiang Univ. Sci. A}}
  \bibinfo{volume}{6} (\bibinfo{year}{2005}), \bibinfo{pages}{116--123}.
\newblock


\bibitem[\protect\citeauthoryear{Speleers}{Speleers}{2019}]%
        {Speleers:2019}
\bibfield{author}{\bibinfo{person}{H. Speleers}.}
  \bibinfo{year}{2019}\natexlab{}.
\newblock \showarticletitle{Algorithm 999: Computation of multi-degree
  {B}-splines}.
\newblock \bibinfo{journal}{{\it ACM Trans. Math. Software}}
  \bibinfo{volume}{45}, Article \bibinfo{articleno}{43} (\bibinfo{year}{2019}),
  \bibinfo{numpages}{15}~pages.
\newblock


\bibitem[\protect\citeauthoryear{Toshniwal, Speleers, Hiemstra, and
  Hughes}{Toshniwal et~al\mbox{.}}{2017}]%
        {Toshniwal:2017polar}
\bibfield{author}{\bibinfo{person}{D. Toshniwal}, \bibinfo{person}{H.
  Speleers}, \bibinfo{person}{R.~R. Hiemstra}, {and} \bibinfo{person}{T.~J.~R.
  Hughes}.} \bibinfo{year}{2017}\natexlab{}.
\newblock \showarticletitle{Multi-degree smooth polar splines: A framework for
  geometric modeling and isogeometric analysis}.
\newblock \bibinfo{journal}{{\it Comput. Methods Appl. Mech. Eng.}}
  \bibinfo{volume}{316} (\bibinfo{year}{2017}), \bibinfo{pages}{1005--1061}.
\newblock


\bibitem[\protect\citeauthoryear{Toshniwal, Speleers, Hiemstra, Manni, and
  R.}{Toshniwal et~al\mbox{.}}{2020}]%
        {Toshniwal:2020}
\bibfield{author}{\bibinfo{person}{D. Toshniwal}, \bibinfo{person}{H.
  Speleers}, \bibinfo{person}{R.~R. Hiemstra}, \bibinfo{person}{C. Manni},
  {and} \bibinfo{person}{Hughes T.~J. R.}} \bibinfo{year}{2020}\natexlab{}.
\newblock \showarticletitle{Multi-degree {B}-splines: Algorithmic computation
  and properties}.
\newblock \bibinfo{journal}{{\it Comput. Aided Geom. Des.}}
  \bibinfo{volume}{76}, Article \bibinfo{articleno}{101792}
  (\bibinfo{year}{2020}), \bibinfo{numpages}{16}~pages.
\newblock


\bibitem[\protect\citeauthoryear{Unser}{Unser}{2005}]%
        {Unser:2005}
\bibfield{author}{\bibinfo{person}{M. Unser}.} \bibinfo{year}{2005}\natexlab{}.
\newblock \showarticletitle{Cardinal exponential splines: Part {II}---Think
  analog, act digital}.
\newblock \bibinfo{journal}{{\it IEEE Trans. Signal Process.}}
  \bibinfo{volume}{53} (\bibinfo{year}{2005}), \bibinfo{pages}{1439--1449}.
\newblock


\bibitem[\protect\citeauthoryear{Unser and Blu}{Unser and Blu}{2005}]%
        {UnserB:2005}
\bibfield{author}{\bibinfo{person}{M. Unser} {and} \bibinfo{person}{T. Blu}.}
  \bibinfo{year}{2005}\natexlab{}.
\newblock \showarticletitle{Cardinal exponential splines: Part {I}---Theory and
  filtering algorithms}.
\newblock \bibinfo{journal}{{\it IEEE Trans. Signal Process.}}
  \bibinfo{volume}{53} (\bibinfo{year}{2005}), \bibinfo{pages}{1425--1438}.
\newblock


\bibitem[\protect\citeauthoryear{Wang and Fang}{Wang and Fang}{2008}]%
        {Wang:2008}
\bibfield{author}{\bibinfo{person}{G. Wang} {and} \bibinfo{person}{M. Fang}.}
  \bibinfo{year}{2008}\natexlab{}.
\newblock \showarticletitle{Unified and extended form of three types of
  splines}.
\newblock \bibinfo{journal}{{\it J. Comput. Appl. Math.}}
  \bibinfo{volume}{216} (\bibinfo{year}{2008}), \bibinfo{pages}{498--508}.
\newblock


\end{thebibliography}

\end{document}